\documentclass[onefignum,onetabnum]{siamart171218}
\usepackage[utf8]{inputenc}
\usepackage{amsfonts}
\usepackage[english]{babel}
\usepackage[utf8]{inputenc}
\usepackage{amsmath}
\usepackage{graphicx}
\usepackage{amssymb}
\usepackage{mathrsfs}
\usepackage[colorinlistoftodos]{todonotes}
\usepackage{mathtools}
\usepackage{bm}
\usepackage{adjustbox}   
\usepackage{algpseudocode}
\usepackage{algorithm}
\usepackage{subcaption}
\usepackage{tabularx,colortbl}
\usepackage{colortbl}
\usepackage{graphicx}
\usepackage[capitalize]{cleveref}

\usepackage{appendix}
\usepackage{soul}

\usepackage[normalem]{ulem}

\newcommand{\dbar}[1]{\overline{{#1}}}
\renewcommand{\hl}[1]{#1}

\newcommand{\vcr}[1]{\bm{#1}}
\newcommand{\mat}[1]{\bm{#1}}
\newcommand{\tsr}[1]{\pmb{\mathcal{#1}}}
\title{Efficient Tensor Completion Algorithms for Highly Oscillatory Operators}
\author{Navjot Singh\thanks{Department of Computer Science, University of Illinois, Urbana-Champaign, IL, USA.\newline Email: 
		\texttt{\{navjot2,solomon2\}@illinois.edu}} \and Edgar Solomonik\footnotemark[1] \and Xiaoye Sherry Li\thanks{{Applied Mathematics and Computational Research Division}, Lawrence Berkeley National Laboratory, Berkeley, CA, USA.\newline Email: 
		\texttt{\{xsli,liuyangzhuan\}@lbl.gov}} \and Yang Liu\footnotemark[2]}

\date{}

\begin{document}

\maketitle
\begin{abstract}

We address the problem of recovering highly oscillatory operators, represented as $n\times n$ matrices with a fixed set of observed entries. Given that these matrices can be well compressed by butterfly matrix decomposition of $L=\mathcal{O}(\log n)$ levels requiring only $\mathcal \mathcal{O} (n \log n)$ degrees of freedom,  we propose a novel reformulation of the butterfly structure as a compact tensor network. Specifically, we reshape the input matrix as an order $2L+2$ dense tensor, and cast its butterfly decomposition as a tensor network consisting of order $L+3$ dense tensors. 
This enables efficient utilization of the existing software infrastructure for dense and sparse tensor computations.
Next, we propose several tensor completion algorithms based on the tensor reformulation of butterfly format, and compare them against algorithms using the quantized tensor train (QTT) format. These algorithms leverage popular completion methods such as alternating least squares, gradient-based optimization, and alternating directions fitting. Moreover, we propose a novel strategy that uses low-rank matrix completion to efficiently generate an initial guess for the proposed algorithms. To demonstrate the efficiency and applicability of our proposed algorithms, we perform three numerical experiments using simulated oscillatory operators in seismic applications. In these experiments, we use $\mathcal O (n \log n)$ observed entries in the input matrix and demonstrate an $\mathcal{O}(n\log^3 n)$ computational cost of the proposed algorithms, leading to speedups of orders of magnitude and significant accuracy improvements for large matrices compared to state-of-the-art low-rank matrix and QTT tensor completion algorithms.

\end{abstract}

\section{Introduction}
Matrix or tensor completion is the task of inferring missing entries in 2D matrices or higher-order tensors from partial observations, which may be corrupted by noise. Completion problems have been widely studied for applications across data science (e.g., recommendation system \cite{dzhenzher2025low} and computer vision \cite{liu2012tensor}), imaging science (e.g., seismic inversion \cite{9785826,Kumar2015efficient}, magnetic resonance imaging \cite{yi2021joint,lingala2011accelerated}, and X-ray tomography \cite{candes2013matrixcompletion}), and high performance computing \cite{hutter2022high}. 

More recent advances involve integrating more sophisticated low-rank structures and smooth data priors into the completion algorithm frameworks

For highly oscillatory operators such as large-scale Radon transforms in X-ray tomography or Green's functions for high-frequency seismic inversion, the aforementioned low-rank matrix and tensor completion algorithms do not approximate the underlying structure of these operators efficiently, leading to increasing rank(s) and high computational cost. Fortunately, butterfly matrices \cite{michielssen_multilevel_1996,li2015butterfly,liu2021butterfly,Yingzhou_2017_IBF,candemyin09} and tensors \cite{kielstra2024linear} are particularly designed for highly oscillatory operators. They have been used to represent integral transforms \cite{Lexing_SFT_2009,bremer2021rapid,Tygert_2010_spherical,Oneil_2010_specialfunction,Haizhao_2018_Phase}, integral operators \cite{michielssen_multilevel_1996,Han_2017_butterflyLUPEC,Liu_2017_HODBF,Sadeed2022VIE,liu2023fast}, PDE operators, \cite{liu2021sparse,claus2023sparse} and neural network parameters \cite{Khoo2020switchnet,dao2019learning,li2022wide,gonon2024make,DEBUT,fan2022adaptablebutterflyacceleratorattentionbased,chen2022pixelated,pmlr-v162-dao22a,dao2020kaleidoscope,liu2023parameter}. Butterfly decompositions are well-suited to compress measurement operators of sizes $n\times n$ in X-ray and seismic applications with $\mathcal{O} (n\log n)$ (via butterfly matrix) or even $\mathcal{O} (n)$ (via butterfly tensor of \cite{kielstra2024linear}) computational and memory complexities. Therefore, butterfly decompositions will be powerful tools if used in matrix or tensor completion of oscillatory operators, as they require significantly lower complexity and less number of observed entries than other existing low-rank completion algorithms. This motivates us to develop completion algorithms in the butterfly format in this current paper.    

In this paper, we develop efficient algorithms to perform completion of an oscillatory matrix via butterfly decomposition. 
Our contribution is three-fold: (1)
First, we reshape the matrix into a higher dimensional tensor by using similar tensorization techniques as quantized tensor train (QTT) decomposition~\cite{khoromskij2011d}. Under such reshaping, the corresponding butterfly matrix decomposition can be converted into a higher-dimensional tensor network, which we call tensor reformulation of butterfly matrices. Note that the result of tensor reformulation is mathematically equivalent to butterfly matrix representation, but makes it convenient to leverage existing tensor network completion software and kernels. This is not to be confused with the butterfly tensor of \cite{kielstra2024linear} which directly constructs a tensor decomposition from a tensor operator without reshaping. (2) Next, we propose three tensor completion algorithms for the tensorized butterfly matrix based on ALS, ADF and gradient optimization, with an initial guess generated and converted from low-rank matrix completion algorithms. We analyze and demonstrate that for highly oscillatory operators, in the experiments presented in this manuscript, these proposed completion algorithms require only $\mathcal{O}(n\log n)$ observed entries and thus achieve $\mathcal{O}(n\log^3 n)$ computational complexities. This is in stark contrast to other completion algorithms such as QTT completion or low-rank matrix completion, 
which require asymptotically larger numbers of arithmetic operations. It is also worth mentioning that although ALS and gradient optimization-based algorithms have been considered for butterfly matrices
\cite{le2025butterflyfactorizationerrorguarantees,DEBUT,dao2019learning}
, these existing algorithms assume $n^2$ observed entries instead of a matrix completion setting.
(3) Finally, we apply the proposed butterfly completion algorithms to several operators in seismic imaging applications to numerically demonstrate their superior convergence, accuracy and efficiency compared with other existing completion algorithms. Aside from applications in tensor completion, we also argue that the proposed tensor reformulation of butterfly can be used to represent neural network layers and allows for efficient automatic differentiation.

The rest of the paper is organized as follows. In \cref{sec:background_tensor} we introduce background information about tensor decomposition, completion and existing software. In \cref{sec:matrix_butterfly} we briefly summarize the definition of the butterfly matrix. The proposed tensor reformulation of a butterfly matrix is detailed in \cref{sec:tensor_butterfly} followed by two proposed completion algorithms with implementation detail and complexity analysis in \cref{sec:butterfly_completion}. Finally, we provide three sets of numerical examples to demonstrate the accuracy, convergence and efficiency of the proposed algorithms in \cref{sec:experiments}.   

%
%

\subsection{Notations}
Given a matrix $\mat{T}$, we use $\mat{T}(i,j)$ or $\mat{T}(t,s)$ to denote its entries or submatrices, where $i,j$ are indices and $t,s$ are index sets. Alternatively, we use $\vcr{t}_i$ to denote its $i$-th row or column, and $t_{i,j}$ to denote its entries. For an order $d$ tensor $\tsr{T}$, we use $\tsr{T}(i_1,i_2,\cdots,i_d)$ to denote its entries, where $(i_1,i_2,\cdots,i_d)$ are index tuples. Alternatively, we use $t_{i_1,i_2,\cdots,i_d}$ to denote its entries. For a set $\Omega$ consisting of order $d$ index tuples, $\Omega_{i_k}$ denotes its subset consisting of tuples whose $k$-th elements are fixed to $i_k$.
Throughout this chapter, for brevity of notation, we use superscripts and subscripts  to label tensors, e.g., $\mat{A}^k$ refers to a matrix with label $k$ and not the $k$th power of $\mat{A}$, $\mat A_{i,j}$ refers to a matrix indexed by two indices $i,j$ and does not refer to its row and columns. 
We define a block diagonal matrix with rectangular (possibly non-square) blocks using the notation $
\text{diag}_i\big\{\mat{T}_{i}\big\}_{1}^{n}$ 
to represent a block diagonal matrix whose blocks are the matrices \( \mat{T}_{i} \in \mathbb{R}^{m_i \times n_i} \), for \( i = 1, \ldots, n-1 \). That is, it is a matrix given as
\begin{equation}
\text{diag}_i\big\{\mat{T}_{i}\big\}_{1}^{n} := 
\begin{bmatrix}
\mat{T}_{1} & & & \\
& \mat{T}_{2} & & \\
& & \ddots & \\
& & & \mat{T}_{n-1}
\end{bmatrix},
\end{equation}
where all off-diagonal blocks are zero. 
The submatrix corresponding to block \( \mat{T}_{i} \) occupies rows from \( \sum_{l=1}^{i-1} m_l +1 \) to \( \sum_{l=1}^{i} m_l  \), and columns from \( \sum_{l=1}^{i-1} n_l +1 \) to \( \sum_{l=1}^{i} n_l  \) with $1-$based indexing.

\section{ Background on Tensor Computations}\label{sec:background_tensor}
\subsection{Tensor Decompositions}

Tensor decompositions are mathematical constructs used to efficiently represent and manipulate high-dimensional tensors. The key idea behind tensor decompositions is to decompose a tensor into a set of simpler tensors linked together according to a specific pattern or topology. A tensor decomposition is represented by a set of parameter tensors $\mathcal S$ that are contracted in a specific way according to the definition of the decomposition to reconstruct the tensor $\tsr X(\mathcal{S})$.

The simplest example of a low rank tensor decomposition is a low rank matrix decomposition of a matrix $\mat M \in \mathbb{R}^{I\times J}$ with rank $R$, given by $\mat M \approx \mat S^{1} \mat S^{2}{}^T$, with parameter
tensors
$\mat S^{1} \in \mathbb{R}^{I \times R}$, $\mat S^{2} \in \mathbb{R}^{J \times R}$, and parameter tensor set $\mathcal{S} = \{ \mat S^{1}, \mat S^{2}\}$, where $\mat X (\mathcal S) = \mat S^{1} \mat S^{2}{}^T$.

Computing an optimal rank tensor decomposition for most problems is NP-hard~\cite{hillar2013most}. Therefore, numerical optimization techniques are used to compute the tensor decomposition by formulating an objective function that is minimized with respect to the tensor decomposition parameters. We formally define the low-rank tensor decomposition optimization problem for a given tensor $\tsr T$ in this subsection. The objective that is minimized to compute the low rank tensor decomposition $\tsr X(\mathcal{S})$,
\begin{align}
\label{eq:gen_fro_dec0}
    \min_{ \mathcal{S} } \ \frac{1}{2} \left\| \tsr T - \tsr X( \mathcal{S}) \right\|_F^2,
\end{align}
where $\tsr X (\mathcal S)$ is the reconstructed tensor. Another critical parameter that is to be determined and not mentioned here is rank(s) or dimensions of the tensors in the set $\mathcal S$. These parameters are often determined based on the desired size and representation accuracy of the tensor decomposition.  

Various types of tensor decompositions have been proposed such as CANDECOMP/PARAFAC (CP)~\cite{harshman1970foundations}, Tucker~\cite{de2000multilinear}, tensor train~\cite{oseledets2011tensor}, hierarchical Tucker~\cite{hackbusch2009new}, tensor ring~\cite{zhao2016tensor}. We review the tensor train decomposition in this paper as it has been shown to approximate oscillatory operators~\cite{dolgov2012superfast,chen2024direct,kielstra2024linear} and is considered as one of the baseline algorithms in this work. Given an order $2L+2$ tensor $\tsr T$ representing a $L+1$ dimensional operator, the tensor train decomposition or matrix product operator (MPO) representation \cite{hubig2017generic} of $\tsr T$ is written as  
\begin{align}
t_{i_0, \dots, i_L, j_0 \dots  j_L} = \sum_{r_1, \dots, r_{L}}s^1_{i_0,j_0,r_1}\Big(\prod_{m = 1}^{L-1
}s^{m +1}_{i_m, j_{m}, r_{m},r_{m+1}} \Big)s^{L+1}_{i_L,j_L,r_L}
\label{eq:QTT_rep}
\end{align}
Here $i_0,\ldots,i_L$ and $j_0,\ldots,j_L$ are used to index the input and output dimensions of the operator $\tsr T$, and the parameter tensor set is $\mathcal{S}=\{\tsr S^1, \tsr S^2,\ldots,\tsr S^{L+1}\}$. Note that the tensors $\tsr S^{1}$ and $\tsr S^{L+1}$ are order $3$ tensors, and the rest are order $4$ tensors. When given an operator in the matrix form, e.g., $\mat T \in \mathbb{R}^{n \times n}$ with typically $n=2^{L+1}$, the matrix can be reshaped and permuted into an order $2L+2$ tensor and then represented by the tensor train decomposition in~\cref{eq:QTT_rep}. This is called the quantized tensor train decomposition (QTT)~\cite{khoromskij2011d} and has been used to represent both sparse and dense matrices~\cite{oseledetsmatrixapproximation,oseledetstensorinv}. Tensor train decomposition (specifically MPO) has been widely used in accelerating quantum chemistry simulations ~\cite{keller2015efficient}, quantum circuit simulations ~\cite{gelss2022low}, and many time-dependent, steady-state, low-dimensional and high-dimensional PDEs~\cite{gelss2022solving,CORONA2017145,manzini2023tensor,ye2024quantized,dolgov2012fast}.

\subsection{Tensor Completion}
The objective function in~\cref{eq:gen_fro_dec0} assumes that all entries of the tensor are available and utilized in computing the decomposition.
However,  in  many applications of interest, only a fixed and limited subset of the tensor entries 
is observed. To address this, a more general objective function is considered, given as
\begin{align}
\label{eq:gen_fro_dec1}
    \min_{  \mathcal S } \frac 1 2\| \mathcal{P}_{\Omega}\Big(\tsr T - \tsr X ( \mathcal S) \Big)\|_F^2,
\end{align}
where $\Omega$ is the set of observed entries, and $\mathcal{P}_{\Omega}$ is a selection operator that selects the entries of tensor denoted by $\Omega$. When $\Omega$ is a strict subset of all the entries of the tensor, the computation of the tensor decomposition and using the tensor decomposition to compute the unobserved entries of the tensor is commonly referred to as \emph{tensor completion}. Low rank matrix completion~\cite{jain2013low,hastie2015matrix,vandereycken2013low} is one of the most established subjects in this category. Besides low rank matrix completion, many algorithms have been proposed for tensor completion in different formats, such as CP~\cite{liu2012tensor}, Tucker~\cite{filipovic2015tucker,kressner2014low}, tensor train~\cite{grasedyck2015variants,cai2022provable,steinlechner2016riemannian}, hierarchical Tucker~\cite{da2015optimization,rauhut2015tensor}, and tensor ring decomposition~\cite{wang2017efficient}.

It is important to note that tensor completion differs from scenarios where arbitrary entries of the tensor can be sampled at will. Instead, only a \textit{fixed} subset of entries are provided. The arbitrary sampling problem is generally easier than the completion setting, and many algorithms have been proposed for different tensor decompositions in this context as well~\cite{oseledets2010tt, saibaba2016hoid, dolgov2021functional}.

\subsection{Software for Tensor Decomposition and Completion}
Computing tensor decomposition and completion involve optimization of objective functions in ~\cref{eq:gen_fro_dec0} and~\cref{eq:gen_fro_dec1}. There is a rich space of algorithms for these optimization problems,
but their primary computational components are tensor contractions and linear algebra.
Hence, the algorithms can be implemented efficiently by making use of the existing software libraries and infrastructure. We review several such libraries below.


Tensor Toolbox~\cite{bader2008efficient} is a MATLAB-based library that provides a wide variety of tensor operations and decompositions for both dense and sparse data. It supports popular models such as CP, Tucker, and block term decompositions. Due to its MATLAB foundation, Tensor Toolbox is widely used for prototyping and algorithm development. However, it is limited to single-node computation and thus is not suitable for very large-scale problems.

Cyclops Tensor Framework (CTF)~\cite{ctf} is a C++ library that provides software infrastructure to perform dense and sparse large-scale tensor  computations in a distributed-memory setting. CTF has been used to implement various tensor decomposition algorithms~\cite{ma2018accelerating, ma2021efficient,singh2021comparison,singh2023alternating}. Recently, new sparse tensor kernels to support tensor completion algorithms in Python were added to CTF~\cite{singh2022distributed}. Additionally, a framework to support sparse tensor contractions with tensor networks~\cite{kanakagiri2024minimum} has been proposed using CTF.

SPLATT (Sparse Tensor Algebra Toolbox)~\cite{smith2015splatt} is a high-performance C/C++ library specifically optimized for sparse tensors. It provides parallel implementations of least squares CP decomposition as well as CP completion algorithms and supports modes of parallelism including shared-memory (OpenMP), distributed-memory (MPI), and hybrid configurations. SPLATT uses a compressed sparse fiber (CSF)~\cite{smith2015tensor} data structure that generalizes compressed sparse row (CSR) format to tensors. 

TuckerMPI is a scalable C++/MPI package designed for large-scale tensor compression based on the Tucker decomposition~\cite{Ballard2019TuckerMPI}. It is optimized for distributed computing environments and allows efficient storage and analysis of large-scale dense tensors and grid datasets, making it suitable for large-scale Tucker decomposition.

TensorLy is a Python library that provides high-level tools and efficient implementations for tensor decomposition, tensor algebra, and tensorized deep learning, with extensive backend support for frameworks such as NumPy, PyTorch, and TensorFlow~\cite{Kossaifi2019TensorLy}. Many algorithms for CP, Tucker, tensor train decomposition are available in Tensorly with different backends, allowing for efficient and easy-to-use tensor decomposition on different hardware.

\section{Background on Butterfly Decomposition}\label{sec:matrix_butterfly}

\subsection{Complementary Low Rank Property}
Oscillatory integral operators arise in many areas of applied mathematics and physics, particularly in wave propagation, quantum mechanics, and harmonic analysis. A key challenge in working with such operators is that, when discretized, they lead to matrices that are dense and typically lack global low-rank structure. As a result, standard low-rank compression techniques fail, and naive discretization followed by direct matrix-vector multiplication incurs a high computational cost of $\mathcal{O}(n^2)$, where $n$ is the total number of degrees of freedom. In high-dimensional problems, this cost is expensive. For example, discretizing a 3D domain with $N$ points per dimension leads to $n = N^3$, making the matrix size $N^3 \times N^3$, and the multiplication cost scales as $\mathcal{O}(N^6)$.

This discretized matrix may be dense and high-rank, however, it exhibits a structured pattern in which certain submatrices are numerically low-rank. The property that characterizes which blocks are low-rank is known as the \emph{complementary low-rank property} (CLR)~\cite{li2015butterfly}. 
This property is mathematically described as follows.


Let $\mat{T}$ be an $m \times n$ matrix. Consider rows to be indexed by a set $\chi$ and columns to be indexed by a set $\xi$, where   $\chi = \{0, \dots, m-1\}$ and columns $\xi = \{0, \dots, n-1\}$, where without loss of generality we assume $m = n$.
For each of the row and column sets, 
complete binary
trees $\mathcal{T}_{\chi}$ and $\mathcal{T}_{\xi}$ are associated respectively, each comprising $L = \mathcal{O}(\log n)$ levels. We extend the complete binary trees such that each node at level $L$ has at most $c$ children which referred to as the leaf nodes.
For simplicity, we assume that all nodes have a leaf size $c$, and that the number of rows and columns satisfies $n = c \cdot 2^L$. However, this definition can be readily extended to arbitrary row, column, and leaf-node sizes. Let $\tau^l$ denote the set of all sets of indices at level $l$ in tree $\mathcal{T}_{\chi }$ and similarly $\nu^l$ denote the set of all sets of indices at level $l$ in tree $\mathcal{T}_{\xi}$.
We have that $\forall i,k \in \{0, \dots, 2^{l }-1 \},~l \in \{0,\dots,L\}$ with $i\neq k$, $\tau_i^l\cap \tau_k^l =\emptyset$, $\tau_i^l=\tau_{2i}^{l-1} \cup \tau_{2i+1}^{l-1}$, and $\tau_0^0=\{0,\ldots, n-1\}$ and analogously for $\nu^l$.

The matrix $\mat T$ satisfies \emph{complementary low-rank property} (CLR) if for any row set or node $\tau^l_i$ at level $l$ in $\mathcal{T}_{\chi}$ and column set or node $\nu^{L-l}_j$ at level $L-l$ in $\mathcal{T}_{\xi}$,  
\begin{equation*}
    \mat{T}(\tau^{l}_i,\nu^{L-l}_j) \text{ admits a low-rank approximation}  
\end{equation*}
where the rank of each submatrix $\mat T(\tau^l_i, \nu^{L-l}_j)$, denoted by $r^{l}_{i,j}$, is bounded by a small constant $r$ (the \emph{butterfly rank})\cite{liu2021butterfly,Yingzhou_2017_IBF,li2015butterfly,Pang2020IDBF}. This special property enables the computation a butterfly decomposition, which can reduce the storage complexity of a matrix and computational complexity of matrix-vector multiplications from $\mathcal{O}(n^2)$ to $\mathcal{O}(n \log n)$. We describe the butterfly decomposition and briefly describe the algorithms to compute it in the next subsection.

\subsection{Butterfly Matrix Decomposition}
Let $\mat T$ be an $n \times n$ input matrix with CLR, and $n= c \cdot 2^{L}$, with $L$ being the number of levels as described above. We assume the number of levels, $L= 2M$, to be even for convenience and simplicity of description. A butterfly matrix decomposition can be described for an arbitrary number of levels. By CLR, for any level $l \in \{0,\dots,L\}$ in $\mathcal{T}_{\chi}$, 
\begin{figure}
	\begin{center}
		\begin{tikzpicture}[scale=3.1]
		\def\myshift{.53}
		\def\r{1/32}
		\def\rr{1/16}
		https://www.overleaf.com/project/5e7066b86260fb0001704fc9 \def\h{1/2}
		\def\q{1/4}
		\def\o{1/8}
		\foreach \x in {0,1,...,15} {
			\fill [gray] (\x*\r,1-\x*\rr) rectangle ++(\r,-\rr);
		};
		\draw (0,0) rectangle (\h,1.0);
		\node[above] at (\q,1){$\mat{S}^1$};
		\tikzset{shift={(\myshift,0)}}
		\foreach \x in {0,1,...,7} {
			\fill [gray] (\x*\r,1-\x*\rr) rectangle ++(\r,-\rr);
			\fill [gray] (\q+\x*\r,1-\x*\rr) rectangle ++(\r,-\rr);
		};
		\draw (0,1) rectangle (\h,\h);
		\node[above] at (\q,1){$\mat{S}^2$};
		\tikzset{shift={(\myshift,0)}}
		\foreach \x in {0,1,...,3} {
			\fill [gray] (\x*\r,1-\x*\rr) rectangle ++(\r,-\rr);
			\fill [gray] (\o+\x*\r,1-\x*\rr) rectangle ++(\r,-\rr);
			\fill [gray] (\q+\x*\r,1-\q-\x*\rr) rectangle ++(\r,-\rr);
			\fill [gray] (\h-\o+\x*\r,1-\q-\x*\rr) rectangle ++(\r,-\rr);
		};
		\draw (0,1) rectangle (\h,\h);
		\node[above] at (\q,1){$\mat{S}^3$};
		\draw (0,1) rectangle (\q,1-\q);
		\draw (\h,\h) rectangle (\q,1-\q);
		\tikzset{shift={(\myshift,0)}}

		\foreach \x in {0,1,...,3} {
			\fill [gray] (\x*\rr,1-\x*\r) rectangle ++(\rr,-\r);
			\fill [gray] (\x*\rr,1-\o-\x*\r) rectangle ++(\rr,-\r);
			\fill [gray] (\q+\x*\rr,1-\q-\x*\r) rectangle ++(\rr,-\r);
			\fill [gray] (\q+\x*\rr,1-\h+\o-\x*\r) rectangle ++(\rr,-\r);
		};
		\node[above] at (\q,1){$\mat{S}^4$};
		\draw (0,1) rectangle (\h,\h);
		\draw (0,1) rectangle (\q,1-\q);
		\draw (\h,\h) rectangle (\q,1-\q);
		\tikzset{shift={(\myshift,0)}}
		\foreach \x in {0,1,...,7} {
			\fill [gray] (\x*\rr,1-\x*\r) rectangle ++(\rr,-\r);
			\fill [gray] (\x*\rr,1-\q-\x*\r) rectangle ++(\rr,-\r);
		};
		\draw (0,1) rectangle (\h,\h);
		\node[above] at (\q,1){$\mat{S}^5$};
		\tikzset{shift={(\myshift,0)}}
		\foreach \x in {0,1,...,15} {
			\fill [gray] (\x*\rr,1-\x*\r) rectangle ++(\rr,-\r);
		};
		\draw (0,1.0) rectangle (1,\h);
		\node[above] at (\h,1){$\mat{S}^6$};
		\end{tikzpicture}
	\end{center}
	\caption{Illustration of a $4$-level butterfly matrix decomposition $\mat{X}=\mat{S}^1\mat{S}^2\mat{S}^3\mat{S}^4\mat{S}^5\mat{S}^6$. \label{fig:BF_4L}}
\end{figure}
\begin{equation}
    \label{eqn:IDtaunu}
    \mat{T}(\tau^l_i,\nu^{L-l}_j) \approx \mat{U}^{L-l}_{i,j} \mat{V}^{l}_{i,j}, \quad \forall i \in \{0, \dots, 2^l - 1\},~ \forall j \in \{0, \dots, 2^{L-l} - 1\},
\end{equation}
where \( \mat{U}^{L-l}_{i,j} \in \mathbb{R}^{|\tau^l_i| \times r^{L-l}_{i,j}} \) and \( \mat{V}^{l}_{i,j} \in \mathbb{R}^{r^{L-l}_{i,j} \times |\nu^{L-l}_j|} \) are low-rank decomposition factors.  
Note that for \( l = 0 \), we have \( i = 0 \) only; similarly, for \( l = L \), we have \( j = 0 \) only. Despite possibly being non-unitary, we refer to them as the column and row basis matrices.   
Each of the ranks $r^{l}_{i,j}$ is bounded by a constant $r$, allowing for decomposition of the input matrix into $\mathcal O(\log n)$ block sparse factors, with each factor having $\mathcal O(n)$ nonzeros. The decomposition of the matrix into these factors is known as the butterfly decomposition. 
We derive the form of butterfly decomposition by starting at the finest levels of both the trees and recursively constructing the bases for intermediate levels by propagating the row and column spaces derived at these boundary levels. 

At level $l=0$ in $\mathcal{T}_{\chi}$,
the row basis matrices $\mat{V}^{0}_{0,j}$ each of size $r^{L}_{0,j} \times c$
are defined by a low-rank factorization of each corresponding submatrix $\mat T(\tau_0^0, \nu_j^L)$, $\forall j \in \{0,\dots, 2^L-1\}$. Here $\tau_0^0$ is the row index set at the root of the tree $\mathcal T_{\chi}$. Note that for levels $l\in \{1,\dots,M\}$, each submatrix $\mat T(\tau^l_i,\nu^{L-l}_j)$,
comprises row indices of its children in $\mathcal T_{\chi}$, i.e., the matrix is concatenated column-wise (and split row-wise). As the matrix is concatenated column-wise, the row basis matrices $\mat{V}^{l}_{i,j}$ are formed by concatenating linear combinations of the row basis matrices corresponding to child nodes at the level $l+1$. The row basis matrices $\mat V^{l}_{i,j}$, with corresponding column and row index sets $\tau^l_i$ and $\nu^{L-l}_j$, are represented in a nested fashion in terms of the children of row index sets, $\nu^{L-l+1}_{2j}$ and $\nu^{L-l+1}_{2j+1}$, given as 
\begin{equation}\label{eqn:nested_basis}
\mat{V}^{l}_{i,j} = 
\mat{W}^{l}_{i,j}\begin{bmatrix}
\mat{V}^{l+1}_{i,2j} & \\
& \mat{V}^{l+1}_{i,2j+1}
\end{bmatrix}, \quad \forall l \in \{1,\dots, M\},
\end{equation}
where $\mat{W}_{i,j}^{l}$
is called the transfer matrix for corresponding set $\tau^{l}_i, \nu^{L-l}_{j}$. By CLR, each  $\mat{W}_{i,j}^{l}$ is of size $r_{i,j}^{L-l }\times (r^{L-l +1}_{i,2j} + r^{L-l +1}_{i,2j+1})$, which is at most $r\times2r$.

Analogously, at level $l=0$ of $\mathcal{T}_{\xi}$, the column basis matrices $\mat{U}^{0}_{i,0}$, of size $c \times r^{0}_{i,0}$ 
is defined by a factorization of each corresponding submatrix $\mat T(\tau_i^L, \nu_0^0)$ $\forall i \in \{0,\dots, 2^L-1\}$.
At levels $l \in \{1, \dots, M\}$, the column basis matrices at level $l$ are concatenation of the linear combination of corresponding basis at level $l+1$. Therefore, column basis, $\mat{U}^l_{i,j}$, with corresponding column and row index sets $\tau^l_i$ and $\nu^{L-l}_j$, are represented in a nested fashion in terms of the children of index sets, $\tau^{l+1}_{2i}$ and $\tau^{l+1}_{2i+1}$, given as  
\begin{equation}\label{eqn:nested_basis_P}
\mat{U}^{l}_{i,j} = 
\begin{bmatrix}
\mat{U}^{l+1}_{2i, j} & \\
& \mat{U}^{l+1}_{2i+1, j}
\end{bmatrix}\mat{P}^{l}_{i,j}, \quad \forall l \in \{1,\dots M\},
\end{equation}  
where   $\mat P^{l}_{i,j} $
is also called the transfer matrix for corresponding set $\tau^{l}_i, \nu^{L-l}_{j}$. Each $\mat P^{l}_{i,j}$ is of size $(r^{l-1}_{2i,j} + r^{l-1}_{2i+1,j}) \times r_{i,j}^{l} $ which is bounded by $2r \times r$ by CLR property. 

At level $M = L/2$ of both the trees (half level), from \cref{eqn:IDtaunu}, we have a low-rank representation of a block $\mat T(\tau^M_i, \nu^M_j)$, given by
\begin{align}
    \mat{T}(\tau^M_i,\nu^{M}_j) \approx \mat{U}^M_{i,j} \mat{V}^{M}_{i,j},
\label{eq:middle_level}
\end{align}
where $\mat U^M_{i,j}$ is the column basis matrix at level $M$ and $\mat V^M_{i,j}$ is the row basis matrix at level $M$ of both the trees for a given $i,j$. The column basis matrix $\mat U^M_{i,j}$ can be computed recursively using \cref{eqn:nested_basis_P} corresponding to the given $\tau^M_i$ and $\nu^M_j$. The column basis for $\tau^M_i$ and $\nu_j^M$ is given as
\begin{align}
    \mat U^M_{i,j} &= \dbar{\mat{U}}_{i}\left(\prod_{l=1}^{M}\dbar{\mat{P}}^{l}_{i,j}\right), \nonumber \\
    \text{where }~\dbar{\mat{U}}_i&=\operatorname{diag}_k \Big\{\mat{U}^0_{k,0}\Big\}_{2^{M}i}^{2^{M}(i+1)}
\label{eqn:butterfly_factors_U},
\end{align}
and the transfer factors $\dbar{\mat{P}}^{l}_{i,j}$ for $l \in \{1,\dots,M\}$ are block diagonal matrices. They consist of transfer matrices $\mat{P}^l_{i,j}$, given as
\begin{align}
\dbar{\mat{P}}^{l}_{i,j}&=\operatorname{diag}_k\Big\{\mat{P}_{k,j}^{l} \Big\}_{2^{M-l}i}^{2^{M-l}(i+1)}.
\label{eqn:butterfly_factors_P}
\end{align}
The row basis matrix $\mat V^M_{i,j}$ can be computed recursively using \cref{eqn:nested_basis} and described analogously.

Since~\cref{eq:middle_level} is defined for an arbitrary $i$ and $j$, collecting all the submatrices of $\mat T$ at level $M$, we have the following approximation of $\mat T$,
\begin{align}
\mat{T} \approx  \mat K = \mat{S}^1\bigg(\prod_{m=1}^{M}\mat{S}^{m+1}\bigg)\bigg(\prod_{m=M+1}^{L}\mat{S}^{m+1}\bigg)\mat{S}^{L+2},\label{eq:BF_mat}
\end{align}
where, 
\begin{align}
\mat{S}^1&=\operatorname{diag}_i\Big\{ \mat U ^{0}_{i,0}\Big\}_{0}^{2^L},\nonumber\\
(\mat{S}^{m+1})^T&=\operatorname{diag}_{j}\Bigg \{\begin{bmatrix}
\operatorname{diag}_{i}\big\{(\mat{P}^{m}_{i,2^{m-1}j})^T\big\}_{0}^{2^{L-m}}\\
\operatorname{diag}_{i}\big\{(\mat{P}^{m}_{i,2^{m-1}j+1})^T\big\}_{0}^{2^{L-m}}  
\end{bmatrix}\Bigg\}_{0}^{2^{m-1} },~\forall m \in \{1,\dots, M \} 
\label{eq:BF_mat_S},
\end{align}
and $\mat S^{L+2}$, and $\mat S^{m+1},~m \in \{M+2,\dots, L\}$ are defined analogously.

 For reference, we plot an example of 4-level butterfly matrix decomposition of \cref{eq:BF_mat_S} in \cref{fig:BF_4L}. In $\mat{S}^2$ and $\mat{S}^3$, each approximately $2r\times r$ block (shown in gray) represents one $\mat{P}^{l}_{i,j}$ in $\mat{S}^4$ and $\mat{S}^5$, each approximately $r\times 2r$ block (shown in gray) represents one $\mat{W}^l_{i,j}$.

Each $\mat{U}^L_{i,0}$ and $(\mat{V}^{0}_{0,j})^T$ is of size $c \times r$ with $i,j \in \{0,\dots,2^L - 1\}$.  
Each $\mat{P}^m_{i,j}$ is of size $2r \times r$ with $i \in \{0, \dots, 2^{L - m} - 1\}$ and $j \in \{0,\dots,2^m - 1\}$ for $m \in \{1,\dots, M\}$.  
Similarly, each $\mat{W}^l_{i,j}$ is of size $r \times 2r$, with $i \in \{0,\dots,2^l - 1\}$ and $j \in \{0, \dots, 2^{L - l} - 1\}$ for $l \in \{1,\dots, H\}$. Therefore, the outermost blocks require $O(cr \cdot 2^L)$ memory, and all of the inner blocks require $O(r^2 L \cdot 2^{L+1})$ memory.  
Assuming $r$ and $c$ to be constant, and $L = O(\log n)$, the full butterfly matrix decomposition can be represented in $O(n \log n)$ memory. We present a brief overview of algorithms to compute the butterfly decomposition below.

\paragraph{Computation of butterfly decomposition} When any entry of $\mat T$ can be sampled in $\mathcal O (1)$, a butterfly decomposition can be efficiently approximated using SVD~\cite{li2015butterfly} in $\mathcal O(n^{1.5})$, or using
Chebyshev interpolation \cite{candes2007fast,Pang2020IDBF} or proxy points \cite{Eric_1994_butterfly,Han_2017_butterflyLUPEC,michielssen_multilevel_1996} in $\mathcal O(n\log n)$ complexity. When $\mat{T}$ can be applied to arbitrary random vectors in $\mathcal O(n \log n)$ complexity, a butterfly decomposition can be approximated using sketching algorithms~\cite{liu2021butterfly} in $\mathcal O(n^{1.5}\log n)$ complexity. Note that none of these methods can be applied when only a fixed subset of entries is given (i.e., the matrix completion setting). In this setting, numerical optimization algorithms must be used and are described in~\cref{sec:butterfly_completion}. 

\subsection{Butterfly Decomposition in Machine Learning}
Due to its ability to approximate Fourier integral operators with a small number of parameters, butterfly decomposition has recently been integrated into deep neural networks to efficiently learn structured linear transforms~\cite{dao2019learning,Fast-transform}. Building on this, the Kaleidoscope hierarchy, which is a product of learnable orthogonal butterfly matrices, was introduced into the Mel-Frequency Spectral Coefficients (MFSC) pipeline~\cite{davis1980comparison}, achieving performance comparable to hand-crafted transforms. These Kaleidoscope matrices are capable of learning latent permutations and accelerating inference~\cite{dao2020kaleidoscope}. To further enhance their practicality and compatibility with GPU architectures, the sparsity structure of the butterfly factors was modified in various works~\cite{DEBUT,pmlr-v162-dao22a,chen2022pixelated}. Empirical evaluations in these works show that restructuring maintains similar accuracy while improving efficiency by enhancing the ease of implementation.

Integrating butterfly matrices into machine learning frameworks like PyTorch has led to new developments in neural network architectures and computational kernels. Although butterfly decompositions theoretically enable 
 matrix-vector products in $O(n\log n)$ complexity, achieving these speedups in practice poses significant challenges. Traditional block-sparse representations often lead to inefficient high-level implementations, as they require manual looping over sparse blocks, which impedes performance. Realizing theoretical speedups for both the forward pass, that is, applying the butterfly matrix, and backward pass, which is obtaining the gradients of block sparse matrices, demands carefully optimized software kernels. While recent efforts have developed efficient GPU kernels for the forward pass~\cite{gonon2024make,fan2022adaptablebutterflyacceleratorattentionbased}, kernels for computing gradients (i.e., the backward pass) remain underdeveloped. To address this limitation, some works have modified the sparsity structure of the butterfly factors to better align with hardware efficiency~\cite{pmlr-v162-dao22a}. In contrast, we propose 
a new representation of the butterfly decomposition as a tensor network, enabling differentiation with respect to the factors by leveraging automatic differentiation with 
tensor contractions. This reformulation is detailed in the next section. 
\begin{figure}
	\centering
		\includegraphics[width=0.9\linewidth]{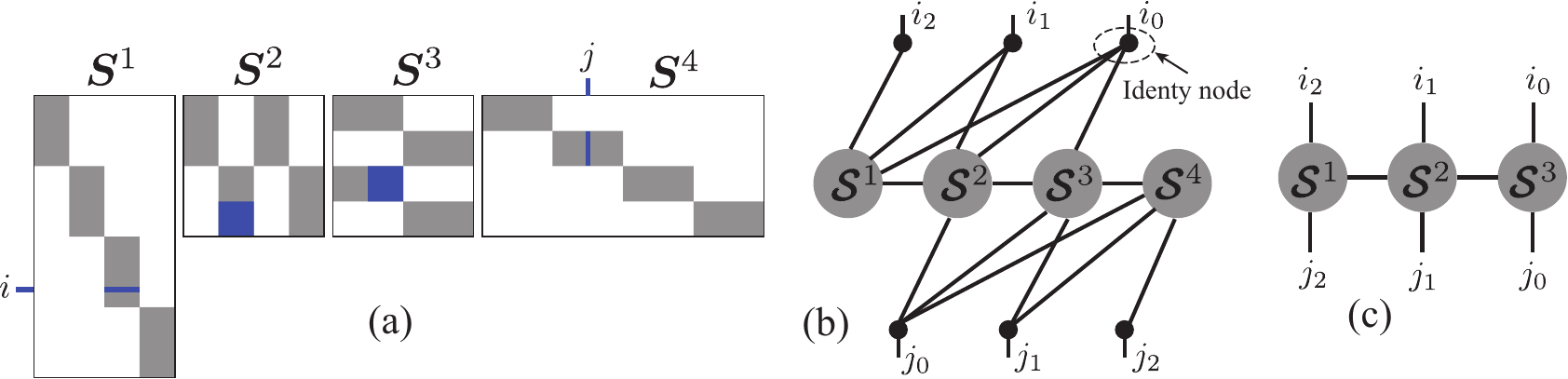}    
	\caption{(a) A $2$-level butterfly matrix decomposition $\mat{X}=\mat{S}^1\mat{S}^2\mat{S}^3\mat{S}^4$. For a matrix entry $(i,j)$ with $i=2c+a$ and $j=c+b$, $a,b\leq c$, $x_{ij}$ is computed by the product of the subblocks highlighted in blue. (b) The diagram of reshaping the $2$-level butterfly matrix decomposition $\mat{X}$ as an order $2L+2$ tensor $\tsr{X}$ with tensor cores $\tsr{S}^1,\tsr{S}^2,\cdots,\tsr{S}^{L+1},\tsr{S}^{L+2}$. Each gray circle represents one tensor core $\tsr{S}^l$ of order $L+3$ for $1<l<L+2$ or order $L+2$ otherwise. Each black node denotes an identity node, representing an external edge $i_0$, $i_1$, $i_2$ or $j_0$, $j_1$, $j_2$. The $(i,j)$ entry in (a) now becomes $(i_0,i_1,i_2,j_0,j_1,j_2)=(1,0,a,0,1,b)$. (c) The diagram of the QTT representation of the same tensor in (b): $\tsr{X}$ with tensor cores $\tsr{S}^1, \tsr{S}^2, \tsr{S}^3$ where each tensor core has order 3 or 4, independent of $L$. \label{fig:tensor_diagram}}
\end{figure}

\section{Tensor Network for Butterfly Decomposition}\label{sec:tensor_butterfly}
In this section, we propose a reformulation that recasts the butterfly decomposition as a tensor decomposition. Each butterfly factor is modeled as a dense tensor, and the product of block sparse matrices in the butterfly decomposition is expressed through tensor contractions. These contractions are naturally supported by modern auto-differentiation frameworks via einsum-like primitives, eliminating the need for hand-crafted software kernels. 

\subsection{Butterfly Decomposition Rewritten as a Tensor Network}
To describe the tensor decomposition reformulation, we will continue to use the same assumptions as in the previous section. We consider an $m \times n$ input matrix $\mat T$, where $n = c \cdot 2^L$, and $L$ is the number of levels in the butterfly decomposition. The row and column sets, $\chi = \{0, \dots, m-1\}$ for rows and $\xi = \{0, \dots, n-1\}$ for columns of the matrix, and $m = n$. These sets are organized into hierarchical trees $\mathcal{T}_{\chi}$ and $\mathcal{T}_{\xi}$ respectively as described in Section~\ref{sec:matrix_butterfly}. 
$\tau^l$ denotes the set of all nodes at level $l$ in tree $\mathcal{T}_{\chi }$ and $\nu^l$ denotes the set of all nodes at level $l$ in tree $\mathcal{T}_{\xi}$. Additionally, we assume that for each level $l \in \{0,\dots, L\}$, the ranks of all the blocks are equal, that is, 
$\forall i,j, \ r^{l}_{i,j}= R_{l}$. For the ease of later complexity analysis, we further assume $R_l=r$.

We reshape the matrix $\mat T$, indexed by $(i,j)$—where $i$ denotes the row and $j$ denotes the column—into a tensor. The tensor is formed by associating an $(L+1)$-tuple to each of $i$ and $j$, based on the nodes $\tau^l_i$ and $\nu^l_j$ that contain $i$ and $j$ respectively, at each level $l$. For generating the tuple corresponding to rows, let $\tau^{0}_0$ be the root node $\rho$. For each level $l$ from $0$ to $L-1$:
\begin{itemize}
    \item Assign $i_l=0$ if $i$ is in the left child of $\rho$ else assign $i_l =1$.
    \item Update the root $\rho$ to be its child where $i$ is found.
\end{itemize}
This recursive process generates a binary tuple $(i_0,  \ldots, i_{L-1})$ for each index $i$. The last index of the tuple, $i_{L}$, is generated by finding the position out of $\{0,\dots,c-1\}$ of $i$ at the leaf level. The reshaping of index $i$ into the tuple $(i_0, \dots, i_L)$ is formally defined by
\begin{align}
    i_l = 
    \begin{cases} 
        i \bmod c, & \text{if } l = L, \\
         \left\lfloor\frac{i }{c \cdot2^{L-l-1}}\right 
    \rfloor  \bmod 2, & \text{if } l = 0, \ldots, L-1.
    \end{cases}
    \label{eq:indices}
\end{align}
A transformation to the index $j \rightarrow (j_0,\dots,j_L)$ is defined analogously. This leads to reshaping of the matrix $\mat T$ into an order $2L + 2$ tensor, $\tsr T$, by representing each entry of the matrix $(i,j)$ by a $(2L +2)-$ tuple  $\underbrace{(i_0, \dots, i_d}_{i}, \underbrace{j_0, \dots, j_d)}_j$.

With this approach, each of the block sparse butterfly factor matrices, $\mat S^{l}$ $\forall l \in \{1, \dots, L+2 \}$, given in~\eqref{eq:BF_mat_S} are transformed into higher-order dense tensors $\tsr S^l$. This is done by stacking up the non-zero blocks of the block-sparse matrices in $\mat S^{l}$ in~\eqref{eq:BF_mat_S} as slices of a higher-order tensor $\tsr S^l$.

To define this transformation, we define a function $\psi$. When given an integer $i \in \{0,\dots,2^{l}-1\}$, it outputs a binary tuple as follows,
\begin{align}
\psi(i,l) = \left( b_0, b_1, \dots, b_{l-1} \right) \in \{0,1\}^l,
\label{eq:function_ind}
\end{align}
where each \( b_m \in \{0,1\} \) is given by
\begin{align}
b_m = \left\lfloor \frac{i}{2^m} \right\rfloor \bmod 2, \quad \text{for } m = 0, 1, \dots, l-1.
\label{eq:fn_ind2}
\end{align}
We can describe the transformation for the butterfly factor matrices $\mat S^{l}$ to tensors $\big\{ \tsr S^l : l \in \{1,\dots,L+2\} \big\}$, using~\eqref{eq:indices},~\eqref{eq:function_ind}, and~\eqref{eq:fn_ind2} with $\tsr S^{1}$ and $\tsr S^{L+2}$ being order $L+2$  and $\big \{\tsr S^{l} : l \in \{2,\dots L+1\} \big\}$ being order $L+3$ tensors. The slices of the higher-order tensors are indexed as follows
\begin{align}
    \forall i, \ \mat S^1 (i_0,\dots,i_{L-1}, :,  :) &= \mat U^0_{i,0}, \ \ (i_0,\dots,i_{L-1}) = \psi(i,L)\nonumber,\\
    \forall i, j, \ \mat S^{m+1}(i_0 \dots i_{L-m} , j_0, \dots, j_L,:,:) &= \mat P^{m}_{i,j},\nonumber \\
    (i_0,\dots,i_{L-m}) = \psi(i,L-m), \ \ (j_0,\dots j_m) &= \psi(j,m), \ \ \forall m \in \{1,\dots,M\},
    \label{eq:tensor_form_factors}
\end{align}
where $\mat U^L_{i,0}$ are the column basis in~\cref{eqn:butterfly_factors_U}, and $\mat P^m_{i,j}$ are column transfer matrices in~\cref{eqn:butterfly_factors_P}. The tensors $\big\{\tsr S^{l} : l \in \{M+2,\dots, L+2\} \big\},$ are defined analogously.

Using the above assignment, each entry of the matrix $\mat T$ in~\cref{eq:BF_mat} is mapped to a tensor $\tsr T$. The input tensor $\tsr T$ is approximated by a butterfly matrix $\mat K$, which is reshaped into $\tsr K$, and is represented by the tensor network with parameters $\mathcal S = \big \{ \tsr S^l : l \in \{1,\dots,L+2\} \big \}$. Using the transformations above, each $c \times c$ block in~\cref{eq:BF_mat} is
\begin{align}
    \mat K(\psi(i,L), :, \psi(j,L), :)& = \mat S^1(\psi_L(i,L),:,:)\big(\prod_{m=1}^{M}\mat S^{m+1}(\psi(i, L-1) , \psi(j,m-1),:,:) \big) \nonumber \\
    &\big(\prod_{l=M+1}^{L}\mat S^{l+1}(\psi_L(i,l-1) ,\psi_L(j,L-l),:,:) \big)  \mat S^{L+2}(\psi_L(j,L),:,:),
\label{eq:tensor_rep_before}
\end{align}
where, for convenience, we slightly abuse notation to write 
\begin{align}
    \mat K (i_0,\dots,i_{L-1}, :, j_0, \dots, j_{L-1},  :) = \mat K(\psi(i,L),:, \psi(j,L),:)
\end{align} and similarly for the parameter tensors. This leads to the following elementwise expression,
\begin{align}
t_{i_0, \dots, i_L, j_0 \dots  j_L} &\approx k_{i_0, \dots, i_L, j_0 \dots  j_L} \nonumber \\
&=  \sum_{r_1, \dots, r_{L+1}} s^{1}_{i_0 \dots i_L r_1} \Big(\prod_{m = 1}^{L
}s^{m +1}_{i_0, \dots, i_{L-m}, j_0, \dots, j_{m-1}, r_{m},r_{m+1}} \Big)  s^{L+2}_{j_0 \dots j_{L-1} r_{L+1}j_L}.
\label{eq:BF_tensor_rep}
\end{align}

As an example, \cref{fig:tensor_diagram}(a) shows a 2-level butterfly matrix decomposition $\mat{X}=\mat{S}^1\mat{S}^2\mat{S}^3\mat{S}^4$ and how an entry, $(i,j)$, of the butterfly matrix is generated as highlighted in blue. \cref{fig:tensor_diagram}(b) shows the tensor diagram of the tensor network representation of the butterfly matrix decomposition in \cref{fig:tensor_diagram}(a). Each gray circle represents the order $L+2$ core tensors $\tsr{S}^1$ and $\tsr{S}^4$, or the order $L+3$ core tensors $\tsr{S}^2$ and $\tsr{S}^3$. Each black dot represents an \textit{identity node}, where all edges connected to it represent the same (external) edge.  

For notational convenience, we rewrite~\cref{eq:BF_tensor_rep} as 
\begin{align}
k_{i_0, \dots, i_L, j_0 \dots  j_L} &=  \sum_{r_1, \dots, r_{L+1}}\Big(\prod_{m = 0}^{L+1
}s^{m +1}_{i_0, \dots, i_{L-m}, j_0, \dots, j_{m-1}, r_{m},r_{m+1}} \Big).
\label{eq:BF_tensor_rep_simple}
\end{align}
Here we define $r_{L+2}=j_L$, $r_{0}=i_L$, and extend the $L+2$ order tensors $\tsr S^1$ and $\tsr S^{L+2}$ to $L+3$ order tensors by defining $i_{-1}=1$ and $j_{-1}=1$.

It is worth mentioning that the reformulation of the butterfly decomposition described in~\cref{eq:BF_mat} as a tensor network representation described in~\cref{eq:BF_tensor_rep} is a valid tensor network only when the ranks of all the blocks at each level are the same, and the size of the leaf nodes is also the same. This is due to the fact that only in this case, we are able to stack the dense sub-blocks in butterfly factors and represent them as tensors as in~\cref{eq:BF_tensor_rep}. Also, this tensor network reformulation is not to be confused with approaches that leverage tensor low-rank compression algorithms (e.g., Tucker decomposition) to construct tensorized butterfly algorithms (see \cite{kielstra2024linear} for an example), which leads to a different tensor diagram compared with \cref{fig:tensor_diagram}(b).

\subsection{Tensor Reformulation of Operations on Butterfly Matrices}
Once the butterfly decomposition is represented as a tensor network, the operations involved with butterfly matrices can be expressed using einsum-like primitives. This enables efficient and portable implementation on heterogeneous architectures by leveraging existing tensor contraction libraries. For example, the $O(n\log n)$ complexity of applying a butterfly matrix $\mat X(\mathcal S)$ to a vector $\vcr v$ is realized by first reshaping the matrix and vector as described in~\eqref{eq:indices}. The matrix-vector product $\vcr u = \mat X(\mathcal{S}) \vcr v$ can then be computed by performing the following dense tensor contractions,
\begin{align}
    u_{i_0\dots i_L} &= \sum_{j_0 \dots j_L}x(\mathcal S)_{i_0 \dots i_L j_0 \dots j_L} v_{j_0 \dots j_L}, \nonumber \\
    &= \sum_{r_1, \dots r_{L+1}}\sum_{j_0\dots j_L}\Big(\prod_{m = 0}^{L+1}s^{m+1}_{i_0 \dots i_{L-m}j_0\dots j_{m-1}r_mr_{m+1}} \Big).
    \label{eq:butterfly-mat-vec}
\end{align}

Some algorithms require gradients of the objective function with respect to the butterfly factors for learning with butterfly matrices.  The reformulation of butterfly decomposition as a tensor decomposition enables easy and efficient computation of the gradients via tensor contractions. We derive the gradients of butterfly factors with respect to an arbitrary objective function below.

Consider an objective function $\phi\big(\mat X (\mathcal{S}) \big)$, 
where $\mat X(\mathcal{S})$ is a butterfly matrix with factors $\mathcal{S}$. Expressions for gradients with respect to the factors can be derived by reshaping the butterfly matrix and other parameters involved as in~\eqref{eq:indices}. The Gradient with respect to the resulting tensorized factor $\tsr S^{1}$ is computed via chain rule and is given by
\begin{align}
    \frac{\partial \phi}{\partial s^{1}_{i_0 \dots i_L r_1}} &= 
\sum_{j_0,\dots,j_L}z_{i_0 \dots i_L j_0 \dots j_L}\sum_{r_2, \dots, r_{L+1}}\Big(\prod_{m = 1}^{L+1}s^{m+1}_{i_0 \dots i_{L-m}j_0\dots j_{m-1}r_mr_{m+1}}\Big) \nonumber,\\
z_{i_0 \dots i_L j_0 \dots j_L} &= \frac {\partial \phi}{\partial x_{i_0 \dots i_L j_0 \dots j_L}}.
\label{eq:grad_butterfly_1}
\end{align}
Similarly, the gradient with respect to $\tsr S^{L+2}$ is
\begin{align}
    \frac{\partial \phi}{\partial s^{L+2}_{j_0 \dots j_{L-1} r_{L+1}j_L}} = 
\sum_{i_0,\dots,i_L}z_{i_0 \dots i_L j_0 \dots j_L}\sum_{r_1, \dots, r_{L}}\Big(\prod_{m = 0}^{L}s^{m+1}_{i_0 \dots i_{L-m}j_0\dots j_{m-1}r_mr_{m+1}}\Big),
\label{eq:grad_butterfly_last}
\end{align}
and gradients with respect to $\tsr S^{l+1},  l \in\{1,\dots,L\}$,
\begin{align}
\frac{\partial \phi}{\partial s^{l+1}_{i_0 \dots i_{L-l}j_0\dots j_{l-1} r_lr_{l+1}}} = 
\!\!\!\!\sum_{\substack{i_{L-l+1}, \dots, i_L\\j_l, \dots, j_L}}\!\!\!\!\!\!z_{i_0 \dots i_L j_0 \dots j_L}\!\!\!\!\!\!\!\!\sum_{\substack{r_1, \dots, r_{l-1} \\ r_{l+2}\dots, r_{L+1}}} \!\!\!\!\Big(\prod_{\substack{m = 0\\m\neq l}}^{L+1}s^{m+1}_{i_0 \dots i_{L-m}j_0\dots j_{m-1}r_mr_{m+1}} \Big).
\label{eq:grad_butterfly_l}
\end{align}
These expressions can be used to perform a backward pass for learning with butterfly matrices in neural networks using dense tensor contractions or to obtain gradients for numerical optimization algorithms for tensor completion in the butterfly matrix format using sparse tensor contractions which is detailed in the next section.

\section{Tensor Completion with Butterfly Decomposition}\label{sec:butterfly_completion}
Many applications require computation of a butterfly decomposition when only a fixed subset of matrix entries or indices is given. This computation is formulated as a tensor completion problem. Let $\mat T \in \mathbb{C}^{n \times n}$, $n = c \cdot 2^L$, be the input matrix, with observed entries set $\bar{\Omega}$, then the objective function is given as,
\begin{align}
    \phi(\mathcal{\bar S}) = \frac 1 2 \big\|  \mathcal{P}_{\Omega}\big(\mat T - \mat X(\mathcal{\bar S}\big) \big \|_F^2,
\end{align}
where $\mat X(\mathcal{\bar S})$ is a butterfly matrix with the butterfly factors $\mathcal{\bar S}$ as a set of block sparse matrices as given in~\eqref{eq:BF_mat_S}. Since the numerical optimization of this objective function would require computations with block sparse matrices which require complicated software kernels, we leverage the reformulation of butterfly matrix as a tensor network introduced in~\cref{sec:tensor_butterfly}. We transform the input matrix $\mat T$ into an order $2L +2$ tensor, $\tsr T \in \mathbb{C}^{2 \times \dots \times c \times 2 \times \dots \times c}$ 
by using the transformation of each index $(i,j) \in \bar{\Omega}$ into a tuple $(i_0, \dots, i_L,j_0,\dots j_L ) \in \Omega$ according to~\eqref{eq:indices}. This transformation results in the corresponding tensor completion problem, with objective,
\begin{align}
    \phi(\mathcal{S}) = \frac 1 2 \big\| \mathcal{P}_{\Omega}\big(\tsr T - \tsr X(\mathcal{S}) \big)\big\|_F^2,\label{eq:BF_objective}
\end{align}
where $\tsr X(\mathcal S)$ is the reconstructed tensor. The reconstructed tensor is described by the tensorized butterfly decomposition  dense tensor parameters, $\mathcal{S}$, as described in~\eqref{eq:BF_tensor_rep}. 

Tensor completion problem is a well-studied problem from numerical optimization perspective with many existing algorithms such as alternating least squares (ALS) \cite{paatero1997223,kolda2009tensor}, adaptive moment estimation (ADAM) optimization \cite{kingma2014adam}, 
alternating direction fitting (ADF)
\cite{grasedyck2015variants}, and Riemannian optimization \cite{kressner2014low}. Next, we present the implementation and analyze the computational complexity of 
the proposed butterfly completion algorithms using ALS, ADAM and ADF in
~\cref{subsec:Num_opt}. Several practical considerations required to implement these algorithms are discussed in~\cref{subsec:Imp_details}

\subsection{Optimization for Tensor Completion of Butterfly Matrices}
\label{subsec:Num_opt}
The proposed ALS, ADAM and ADF optimization methods for butterfly completion are described in this subsection. We focus on the explanation of the ALS algorithm and its complexity analysis, and briefly summarize the ADAM and ADF algorithms as many computational kernels are shared between them. 

\subsubsection{Alternating Least Squares (ALS)}
ALS has proven to be a competitive method for tensor completion problems for various tensor decompositions~\cite{kolda2009tensor}. The proposed ALS-based butterfly completion algorithm proceeds by minimizing one factor at a time, summarized in~\cref{alg:ALS_BF}. For minimizing factors of the tensorized butterfly, the resulting least squares  to solve for the first parameter tensor $\tsr S^1$ is given for each $({i_0, \dots, i_{L}})$ and $(j_0, \dots, j_L) \in \Omega_{i_0, \dots, i_{L}}$, as 
\begin{align}
    \sum_{r_1}s^1_{i_0 \dots i_Lr_1} a_{i_0\dots i_{L-1}j_0\dots j_L r_1} =
    t_{i_0 \dots i_L j_0 \dots j_L}, \nonumber\\
    a_{i_0\dots i_{L-1}j_0\dots j_L r_1} = \sum_{r_2, \dots, r_{L+1}}(\prod_{m = 1}^{L+1}s^{m+1}_{i_0 \dots i_{L-m}j_0\dots j_{m-1}r_mr_{m+1}}), 
\end{align}
Similarly, the least squares  to solve for the last parameter tensor, $\tsr S^{L+2}$, is given for each $(j_0, \dots, j_{L})$ and $(i_0, \dots, i_L) \in \Omega_{j_0, \dots, j_{L}}$, as
\begin{align}
    \sum_{r_{L+1}}s^{L+2}_{j_0 \dots j_Lr_{L+1}} a_{i_0\dots i_L j_0 \dots j_{L-1}r_{L+1}} = t_{i_0 \dots i_L j_0 \dots j_L}, \nonumber \\
    a_{i_0\dots i_L j_0 \dots j_{L-1}r_{L+1}} =\sum_{r_1, \dots, r_{L}}(\prod_{m = 0}^{L}s^{m+1}_{i_0 \dots i_{L-m}j_0\dots j_{m-1}r_mr_{m+1}}).
\end{align}
The least squares  to solve for the remaining parameters $\tsr S^{l+1}$, $ l \in \{1,\dots,L\}$, for each $(i_0, \dots, i_{L-l}, j_0, \dots, j_{l-1})$ and $(i_{L-l+1}, \dots, i_L, j_l, \dots, j_L) \in \Omega_{i_0, \dots, i_{L-l}, j_0, \dots, j_{l-1}}$, is given as
\begin{align}
    \sum_{r_l,r_{l+1}}s^{l+1}_{i_0 \dots i_{L-l}j_0 \dots j_{l-1}r_l r_{l+1}} a_{i_0 \dots i_L j_0 \dots j_L r_{l} r_{l+1}} = t_{i_0 \dots i_L j_0 \dots j_L}, \nonumber \\
    a_{i_0 \dots i_L j_0 \dots j_L r_l r_{l+1}} = \sum_{\substack{r_1, \dots, r_{l-1} \\ r_{l+2}\dots, r_{L+1}}}(\prod_{\substack{m = 0\\m\neq l}}^{L+1}s^{m+1}_{i_0 \dots i_{L-m}j_0\dots j_{m-1}r_mr_{m+1}}).
\end{align}
All of the above least squares problems can be solved by solving for fibers of $\tsr S^{(1)}$ and $\tsr S^{L+2}$, and slices of $\tsr S^{l}, l \in \{2,\dots,L+1\}$ in parallel, similar to tensor train completion~\cite{grasedyck2015variants}. We solve the least squares equations by using normal equations that result in small symmetric positive semi-definite systems of size $r \times r$ for $\tsr S^{1}$ and $\tsr S^{L+2}$ or $r^2 \times r^2$ for $\tsr S^{l}, \forall l \in \{2,\dots,L+1\}$, that can be solved in parallel efficiently~\cite{singh2022distributed,smith2016exploration}. We summarize the normal equations method in~\cref{alg:ALS_BF}. This algorithm proceeds by iterating over each fiber of $\tsr S^1$ and $\tsr S^{L+2}$ and slice of $\tsr S^l, l \in \{2,\dots,L+1\}$. This corresponds to lines \ref{line:als:U}, \ref{line:als:V}, and \ref{line:als:ker} and this loop can be parallelized as there is no dependency. These steps and their complexities are explained as follows.

\begin{algorithm}
	\caption{\small \textsf{ALS for Butterfly Completion}\label{alg:ALS_BF}}
    \footnotesize
	\begin{algorithmic}[1]
		\State \textbf{Input:} Observed matrix $\mat{T}_\Omega$ (or its tensorization), rank \( r \), initial guess $\tsr S^{i},i=1,\ldots,L+2$, max iterations \( t_{\max} \), and convergence threshold $\epsilon$.
		\State \textbf{Output:} Updated butterfly cores $\tsr S^{i},i=1,\ldots,L+2$.
		\For{\( t = 1 \) to \( t_{\max} \) and not converged }
		\For {$l = 1$ to $L + 2$ }				
		\If {$l=0$} 
		\For {each $(i_0, \dots, i_{L})$} \Comment{one row of a frontal slice of $\tsr{S}^{1}$}\label{line:als:V}
		\State Initialize $\vcr y = \vcr 0$, $\mat K = \mat 0$
		\For {each $(j_0, \dots, j_L) \in \Omega_{i_0, \dots, i_{L}}$}
		\State $\vcr{v} = 
		\Big(\prod_{m=1}^{L} 
		\mat{S}^{m+1}(i_0, \dots, i_{L-m}, j_0, \dots, j_{m-1}, :,:)
		\Big)  
		\vcr{s}^{L+2}(j_0, \dots, j_L,:)$\label{line:als:v_L2}
		\State Update $\vcr y \mathrel{+}= t_{i_0 \dots i_L j_0 \dots j_L} \vcr{v}$\label{line:als:rhs_L2} 
		\State Update $\mat K \mathrel{+}= \vcr{v}\vcr{v}^T$\label{line:als:K_L2}
		\EndFor
            \State Solve $\mat K \vcr s = \vcr y$  and update $\vcr s^{1}(i_0, \dots, i_{L},:)$  $\leftarrow$ $\vcr s$\label{line:als:solve_L2}
		\EndFor						
		\ElsIf{$l=L+2$} 
		\For {each $(j_0, \dots, j_{L})$} \Comment{one column of a frontal slice of $\tsr{S}^{L+2}$}\label{line:als:U}
		\State Initialize $\vcr y = \vcr 0$, $\mat K = \mat 0$
		\For {each $(i_0, \dots, i_L) \in \Omega_{j_0, \dots, j_{L}}$}
		\State $\vcr{v} = 
		\vcr{s}^{1}(i_0, \dots, i_L,:){}^T\Big(\prod_{m=1}^{L} 
		\mat{S}^{m+1}(i_0, \dots, i_{L-m}, j_0, \dots, j_{m-1},:,:) 
		\Big)$\label{line:als:v_1}
		\State Update $\vcr y \mathrel{+}= t_{i_0 \dots i_L j_0 \dots j_L} \vcr{v}$\label{line:als:rhs_1}  
		\State Update $\mat K \mathrel{+}= \vcr{v}\vcr{v}^T$\label{line:als:K_1}
		\EndFor
            \State Solve $\mat K \vcr s = \vcr y$ and update $\vcr s^{L+2}(j_0, \dots, j_{L},:)$ $\leftarrow$ $\vcr s$ \label{line:als:solve_1}
		\EndFor			
		\Else
		\For {each $(i_0, \dots, i_{L-l}, j_0, \dots, j_{l-1})$} \Comment{one frontal slice of $\tsr{S}^{l+1}$}\label{line:als:ker}
		\State Initialize $\vcr y = \vcr 0$, $\mat K = \mat 0$
		\For {each $(i_{L-l+1}, \dots, i_L, j_l, \dots, j_L) \in \Omega_{i_0, \dots, i_{L-l}, j_0, \dots, j_{l-1}}$}
		\State $\vcr{u} = \vcr{s}^{1}(i_0, \dots, i_L,:){}^T
		\Big(\prod_{m=1}^{l-1} 
		\mat{S}^{m+1}(i_0, \dots, i_{L-m}, j_0, \dots, j_{m-1},:,:)
		\Big) $ \label{line:als:v_l}
		\State $\vcr{v} = 
		\Big(\prod_{m=l+1}^{L} 
		\mat{S}^{m+1}(i_0, \dots, i_{L-m}, j_0, \dots, j_{m-1},:,:)
		\Big)  
		\vcr{s}^{L+2}(j_0, \dots, j_L,:)$\label{line:als:u_l}
		\State Update $\vcr y \mathrel{+}= t_{i_0 \dots i_L j_0 \dots j_L} (\vcr{v}\otimes\vcr{u})$\label{line:als:rhs_l}
		\State Update $\mat K \mathrel{+}= (\vcr{v}\otimes\vcr{u})(\vcr{v}\otimes\vcr{u})^T$\label{line:als:K_l}
		\EndFor
		\State Solve $\mat K \vcr s= \vcr y$  
        and update $\mat S^{l+1}(i_0, \dots, i_{L-l}, j_0, \dots, j_{l-1},:,:) := \mathrm{reshape}(\vcr s,r,r)$\label{line:als:solve_l}
		\EndFor		
		\EndIf
		\EndFor
	\EndFor
	\end{algorithmic}
\end{algorithm}

When solving for $\tsr S^1$ and $\tsr S^{L+2}$, the algorithm then iterates over each observed entry in $\Omega$ and computes the corresponding fibers of the system tensor $\tsr A$ in a vector $\vcr v$ of size $r$ by computing $L$ matrix-vector products with matrices of size $r \times r$ and vectors of size $r$, computed in lines \ref{line:als:v_L2} and \ref{line:als:v_1} with cost $\mathcal O(Lr^2)$. The right-hand side for each fiber is formed in vector $\vcr y$ by accumulating scalar products of $\vcr v$ with corresponding input tensor entries $t_{i_0\dots i_L j_0 \dots j_L}$ in lines \ref{line:als:rhs_L2} and \ref{line:als:rhs_1} with cost $\mathcal{O}(r)$. The semi-definite systems $\mat K$ are computed by accumulating the outer product of vectors $\vcr v$ with themselves in lines \ref{line:als:K_L2} and \ref{line:als:K_1} with cost $\mathcal O (r^2)$. Lines \ref{line:als:v_L2}, \ref{line:als:rhs_L2}, \ref{line:als:K_L2}, \ref{line:als:v_1}, \ref{line:als:rhs_1}, and \ref{line:als:K_1} are executed for each nonzero; therefore, the cost of these operations is $\mathcal  O(|\Omega|Lr^2)$. A positive semi-definite (PSD) solve is performed in lines \ref{line:als:solve_L2} and \ref{line:als:solve_1} to form the updated fiber with cost $\mathcal O (r^3)$. In practice, a regularization parameter is used to ensure the positive definiteness of the system.
Since lines \ref{line:als:solve_L2} and \ref{line:als:solve_1} are executed for all the fibers, the total cost is $\mathcal{O}(nr^3)$. 

When solving for $\tsr S^l, l \in \{2,\dots,L+1\}$, the algorithm then iterates over each observed entry in $\Omega$ and computes the corresponding fibers of the system tensor $\tsr A$ that are a Kronecker product of two vectors $\vcr v$ and $\vcr u$ formed in lines \ref{line:als:v_l} and \ref{line:als:u_l} by performing at most $L$ matrix-vector products, each with cost $\mathcal{O}(r^2)$. Similar to the above case, right-hand sides in line \ref{line:als:rhs_l} and left-hand sides in line \ref{line:als:K_l} are formed by accumulating scalar and outer products of vectors, but of size $\mathcal{O}(r^2)$; each with a cost $\mathcal O(r^2)$ and $\mathcal{O}(r^4)$ respectively. Each slice is then updated by performing a PSD solve in line \ref{line:als:solve_l}, each with a cost $\mathcal O (r^6)$ and reshaping the result into a matrix in row-major ordering. The resulting total cost for updating all the slices is $\mathcal O \big(|\Omega| (Lr^2 + r^4) + nr^6 \big)$. Since the solves are performed for each factor, substituting $L = \mathcal O ( \log n)$, the total cost for one iteration of ALS is $ \mathcal O\big(|\Omega| (r^2\log^2n  + r^4\log n) + nr^6 \log n\big)$. 

The ALS algorithm runs until maximum number of iterations are reached or convergence is achieved. The convergence can be defined by (1) measuring the norm of the difference of factor tensors between iterations, (2) computing the gradient of objective function $\phi(\mathcal{S})$ in  \eqref{eq:BF_objective} or (3) comparing the relative difference of the reconstruction and ground truth using a tolerance $\epsilon$ defined as
\begin{align}
    \big\| \mathcal{P}_{\Omega}\big(\tsr T - \tsr X(\mathcal{S}) \big)\big\|_F/\big\| \mathcal{P}_{\Omega}\big(\tsr T \big)\big\|_F <\epsilon.\label{eq:BF_convergence}
\end{align}
Unless otherwise stated, we use  \eqref{eq:BF_convergence} as the convergence metric in all the completion algorithms under consideration.

Although not detailed here, we remark that many software kernels of~\cref{alg:ALS_BF} can be reused to develop a QTT completion algorithm. The QTT representation of the same $\mat T$ matrix (and its tensorization $\tsr T$) as  \eqref{eq:BF_tensor_rep_simple} is defined as 
\begin{align}
t_{i_0, \dots, i_L, j_0 \dots  j_L} &= \sum_{r_1, \dots, r_{L}}\Big(\prod_{m = 0}^{L
}s^{m +1}_{i_m, j_{m}, r_{m},r_{m+1}} \Big).
\label{eq:QTT_rep_simple}
\end{align}
Comparing with  \eqref{eq:BF_tensor_rep_simple}, QTT can be viewed as a more aggressive compression format than butterfly, consisting of order $4$ tensor cores instead of butterfly's order $L+3$ tensor cores. The diagrams of butterfly and QTT of the same tensor with $L=3$ are shown in \ref{fig:tensor_diagram}(b) and (c). The QTT completion algorithm (based on extension of the tensor train completion in~\cite{grasedyck2015variants}) will be used as one of our reference algorithms in numerical experiments of \cref{sec:experiments} and is discussed in~\cref{alg:ALS_QTT}. We also remark that~\cref{alg:ALS_BF} reduces to low-rank matrix completion algorithm when $L=0$, which is used as one of our reference algorithms in numerical experiments, as well as the initial guess to the butterfly completion algorithms (see \cref{subsec:Imp_details}).

\subsubsection{ADAM Optimization}
Adaptive moment estimation (ADAM)~\cite{kingma2014adam} is a first-order optimization algorithm that combines momentum and adaptive learning rates, making it effective for training deep models and non-convex optimization tasks. ADAM is used in tensor completion problems to efficiently update the factors by obtaining gradients via sparse tensor contractions. The proposed ADAM-based butterfly completion algorithm is described in \cref{alg:ADAM_BF}. For the factors in tensorized butterfly, the gradients are computed by computing the sparse tensor contractions in s~\eqref{eq:grad_butterfly_1},~\eqref{eq:grad_butterfly_l}, and~\eqref{eq:grad_butterfly_last}, with $\tsr Z = \mathcal P_{\Omega} \big(\tsr X (\mathcal S) - \tsr T \big)$. The cost of computing the gradients is $O(|\Omega|Lr^2)$ and a detailed description of the sparse tensor contractions for gradients is summarized in~\cref{alg:ADAM_BF}.  Once gradients are computed, updates can be computed via the ADAM optimization method as summarized in~\cref{alg:Adam_update}.  ADAM uses scalar operations on parameters with the total cost of single iteration being $O( |\Omega|L^2r^2 + Lnr^2)$. Replacing $L= O(\log n)$ in the previous expression, the total cost of a single ADAM iteration is $O(|\Omega|r^2 \log^2n + nr^2\log n)$. 
Comparing with the total cost of a single ALS iteration, ADAM has an $\mathcal{O}(r^2)$ reduction for the leading order term. That said, ADAM oftentimes requires a much larger number of iterations to converge compared with ALS. Just like ALS, we also implement a QTT completion algorithm with ADAM optimization with the same computational kernels as used for butterfly completion.

\begin{algorithm}
\caption{ADAM\_Update($\tsr{S}$,$\tsr{G}$,$\tsr{M}$,$\tsr{V}$,$\alpha$,$\beta_1$,$\beta_2$,$\epsilon$,$t$)}
    \label{alg:Adam_update}
	\begin{algorithmic}[1]
		\State \textbf{Input:} A tensor (or array) of parameters $\tsr{S}$, its gradient $\tsr{G}$, and the current first and second moment estimates $\tsr{M}$ and $\tsr{V}$, learning rate $\alpha$, decay rates for moments $\beta_1,\beta_2$, small constant $\sigma$, current iteration $t$. 
        \State \textbf{Output:} Updated parameters $\tsr{S}$, first and second moment estimates $\tsr{M}$ and $\tsr{V}$			
		\For {each entry $(s,g,m,v)$ in $(\tsr S, \tsr G, \tsr M, \tsr V)$} 
		\State $m \gets \beta_1 m + (1 - \beta_1) g$ \Comment{Update biased first moment estimate}
		\State $v \gets \beta_2 v + (1 - \beta_2) (g \ast g)$ \Comment{Update biased second moment estimate}
		\State $		\hat{m} \gets \frac{m}{1 - \beta_1^t}, \quad 
		\hat{v} \gets \frac{v}{1 - \beta_2^t}$ \Comment{Update bias-corrected estimates}			
		\State $s\gets s - \frac{\alpha \hat{m}}{\sqrt{\hat{v}} + \sigma}$
		\EndFor		
	\end{algorithmic}
\end{algorithm}

\begin{algorithm}
	\caption{\small \textsf{ADAM for Butterfly Completion}\label{alg:ADAM_BF}}
    \footnotesize
	\begin{algorithmic}[1]
		\State \textbf{Input:} Observed matrix $\mat{T}_\Omega$ (or its tensorization), rank \( r \), initial guess $\tsr S^{i},i=1,\ldots,L+2$,, max iterations \( t_{\max} \), learning rate $\alpha$, decay rates for moments $\beta_1,\beta_2$, small constant $\sigma$, and convergence threshold $\epsilon$.
        \State \textbf{Output:} Updated butterfly cores $\tsr S^{i},i=1,\ldots,L+2$.
		\State Initialize gradient tensors, first and second moment estimate tensors, $\tsr G^l =0 \in \mathcal{G}$, $\tsr M^{l}=0 \in \mathcal{M}$, and $\tsr V^{l}=0 \in \mathcal{V}$, $l=1,\ldots,L+2$
		\For{\( t = 1 \) to \( t_{\max} \) and not converged}		
		\For {$l = 1$ to $L + 2$}
            \State Initialize $\tsr G^l =0$
		\If {$l=0$} 
		\For {each $(i_0, \dots, i_{L}) \in \Omega$} \Comment{one row of a frontal slice of $\tsr{S}^{1}$}

                \For {each $(j_0,\dots, j_L) \in \Omega_{i_0\dots i_L}$}
		
		          \State Compute $\vcr g^1(i_0, \dots, i_L,:) +=(x_{i_0 \dots i_L j_0 \dots j_L}- t_{i_0 \dots i_L j_0 \dots j_L}) \vcr{v}$
		          \State where $\vcr{v} = 
		\Bigg(\prod_{m=1}^{L} 
		\mat{S}^{m+1}(i_0, \dots, i_{L-m}, j_0, \dots, j_{m-1},:,:)
		\Bigg)  \vcr{s}^{L+2}(j_0, \dots, j_L,:)$
		      \EndFor
		\EndFor				
		\ElsIf{$l=L+2$} 
		\For {each $(j_0, \dots, j_{L}) \in \Omega $} \Comment{one column of a frontal slice of $\tsr{S}^{L+2}$}
                \For {each $(i_0, \dots, i_L) \in \Omega_{j_0 \dots j_L}$}
		          \State Compute $\vcr{g}^{L+2}(j_0, \dots, j_L,:) +=(x_{i_0 \dots i_L j_0 \dots j_L}- t_{i_0 \dots i_L j_0 \dots j_L}) \vcr{u}$
		          \State where $\vcr{u} = 
		\vcr{s}^{1}(i_0, \dots, i_L,:){}^T\Bigg(\prod_{m=1}^{L} 
		\mat{S}^{m+1}(i_0, \dots, i_{L-m}, j_0, \dots, j_{m-1},:,:)
		\Bigg)$
		      \EndFor
            \EndFor
	\Else
		\For {each $(i_0, \dots, i_{L-l}, j_0, \dots, j_{l-1}) \in \Omega$} \Comment{one frontal slice of $\tsr{S}^{l+1}$}
            \For {each $(i_{L-l+1},\dots i_L,j_{l},\dots,j_L) \in \Omega_{i_0 \dots i_{L-l} j_0 \dots j_{l-1}}$}
		      \State Compute $\vcr{g}^{l+1}(i_0,\dots, i_{L-l},j_0,\dots, j_{l-1},:,:) += (  x_{i_0 \dots i_L j_0 \dots j_L}- t_{i_0 \dots i_L j_0 \dots j_L})  \vcr{v} \otimes \vcr{u}$
		      \State where $\vcr{u} = \vcr{s}^{1}(i_0, \dots, i_L,:){}^T
		\Bigg(\prod_{m=1}^{l-1} 
		\mat{S}^{m+1}(i_0, \dots, i_{L-m}, j_0, \dots, j_{m-1},:,:) 
		\Bigg) $ \\
		      \State $\vcr{v} = 
		\Bigg(\prod_{m=l+1}^{L} 
		\mat{S}^{m+1}(i_0, \dots, i_{L-m}, j_0, \dots, j_{m-1},:,:)
		\Bigg)  
		\vcr{s}^{L+2}(j_0, \dots, j_L,:)$
		      \EndFor	
            \EndFor
		\EndIf
		\State $(\mathcal{S}, \mathcal{M}, \mathcal{V})\leftarrow$ ADAM\_Update($\mathcal{S}$,$\mathcal{G}$,$\mathcal M$, $\mathcal{V}$,$\alpha$,$\beta_1$,$\beta_2$,$\epsilon$,$t$)	
		\EndFor		
		\EndFor
	\end{algorithmic}
\end{algorithm}

\subsubsection{Alternating Direction Fitting (ADF)}
Gradient-based methods like ADAM are computationally cheaper than ALS, however, these methods converge more slowly, taking much more time to reach an accurate decomposition. An alternative to solving least squares systems that arise in the ALS method described above is the ADF algorithm
~\cite{grasedyck2013alternating}. 
ADF can be viewed as a variant of ALS that aims at reducing the cost per iteration while maintaining a similar convergence property to get the best of both worlds
~\cite{grasedyck2015variants}.
The ADF method relies on each of the linear systems being solved for ALS subiteration having orthonormal columns. Once this property is satisfied, ADF replaces these least-squares solves with gradient-based updates of individual core entries with accelerated convergence using an over-relaxation step size computation scheme. Our tensor network reformulation of the butterfly decomposition allows us to make the observation that this property can be satisfied for the butterfly decomposition, enabling a straightforward implementation of the ADF algorithm.

The step-size for each slice of the cores is determined by using the formula described in Algorithm $5$ in
~\cite{grasedyck2015variants}.
The key comparison point with respect to ALS is
that ALS forms normal equations in Line~\ref{line:als:K_l} 
(and Lines
\ref{line:als:K_L2} and \ref{line:als:K_1}) of \Cref{alg:ALS_BF}, 
leading to linear systems of sizes $r^2\times r^2$ requiring $\mathcal{O}(r^4)$ cost for each slice. Instead, ADF  does not compute these systems at all and works with gradient and sparse butterfly reconstruction contractions that cost $\mathcal{O}(|\Omega|r^2 \log^2n)$ per iteration with an additional cost of orthogonalizing cores at each subiteration, which requires $\mathcal{O}(n)$ parallel QR decompositions of matrices of size $2r \times r$. The total per iteration cost of ADF for butterfly thus becomes $\mathcal{O}(|\Omega|r^2 \log^2n + n r^2\log n)$. 
The implementation details of the orthogonalization for ADF for butterfly completion are described in
 ~\cref{app:butterfly-orth}. 
The proposed ADF-based butterfly completion algorithm is described as
\cref{alg:butterfly-adf}.

\subsection{Initial Guess Generation}
\label{subsec:Imp_details}
The ALS, ADAM and ADF methods
require an initial guess for the factors. Although random initial guesses are used for many existing tensor and matrix completion algorithms for data science and scientific computing applications, we remark that they do not work well for the proposed butterfly completion algorithms (
and the reference QTT completion algorithms
), particularly for applications involving highly oscillatory matrices. This is particularly because the larger number of parameters leading to a more complex optimization landscape. Additionally, as we will demonstrate in numerical experiments, we use $|\Omega|=\mathcal{O}(nr \log n)$ observed entries, making full recovery significantly more challenging and requiring an initial guess that is closer to the true solution in Frobenius norm than a random guess. We propose to first use the low-rank matrix completion algorithm (i.e.,~\cref{alg:ALS_BF} with $L=0$ and $r=R$, $R$ being a constant much smaller than the actual rank of the matrix). Neither 
$|\Omega|=\mathcal{O}(nr \log n)$
nor $r=R$ is sufficient for the low-rank completion algorithm to recover the underlying matrix with sufficient accuracy, that is, the relative error over observed entries defined in~\eqref{eq:BF_convergence} stays over $0.5$. However, the resultant low-rank matrix provides a rough approximation of the target matrix and can serve as a better initial guess than a random guess for the butterfly completion algorithm. Our experimental observations suggest that obtaining a low rank decomposition with relative test error of about $0.9$ is enough to finally obtain a good butterfly matrix approximation. It is worth mentioning that the low-rank completion algorithm only requires $\mathcal{O}(t_{\max}|\Omega|R^2 + nR^3)$ operations. 

Once the initial low-rank matrix is computed, the low-rank matrix of rank $R$ needs to be converted to a butterfly matrix of $L$ levels and rank $r$, i.e., $\tsr S^{i},i=1,\ldots,L+2$, to use it as the initial guess. This low-rank to butterfly conversion can be performed efficiently using a randomized algorithm based on a randomized algorithm that uses matrix-vector products of the butterfly matrix~\cite{liu2021butterfly}. In what follows, we briefly describe the randomized algorithm shown in~\cref{alg:low_rank_butterfly_rand}.The input to~\cref{alg:low_rank_butterfly_rand} is the low rank matrix $\mat X = \mat A \mat B$ of rank $R$, desired numbers of butterfly levels $L$ and rank $r$, and an oversampling integer $p$. The algorithm sketches the relevant blocks of $\mat X$ using Kronecker-structured sketching matrices. Specifically, it employs Gaussian random matrices in Kronecker product with appropriately sized identity matrices. We express these operations in tensor form in~\cref{alg:low_rank_butterfly_rand}.

The algorithm proceeds by first reshaping both the low rank factors to tensors using the reshape presented in~\eqref{eq:indices}. We then generate a random Gaussian tensor of appropriate size. The tensor is a reshaping of matrix of size $n \times (r+p)$ using ~\eqref{eq:indices}. A tensor contraction with cost $\mathcal{O}(nrR)$ is computed in Line \ref{line:lr2bf:sketching_outer} to sketch the factors using the generated random Gaussian tensor. The column and row basis for outer-most levels can then be computed with QR factorization with column pivoting over the blocks of the resulting tensors denoted by QRCP in Line \ref{line:lr2bf:qrcp_outer}. Note that these blocks can be executed in parallel, and the QR factorization costs $\mathcal O(nr^2)$ operations, and the blocks of factor tensors are updated to the upper triangular factors obtained from the QR factorization. Next, for each inner level $l=1,\ldots,H$, a new Gaussian random tensor of order $L+2-l$ is generated. The tensors implicitly represent Gaussian random matrices obtained as the Kronecker product of identity matrices and appropriately sized Gaussian random matrices. The sketching and QR factorization are performed at Line \ref{line:lr2bf:sketching_inner}, \ref{line:lr2bf:qrcp_inner}. These operations respectively cost $\mathcal{O}(nRr^2)$, $\mathcal O (nr^3)$ for each level $l$. Once all factor cores have been computed, $\tsr S^{H+1}$ is updated by contracting the resulting tensors in Line \ref{line:lr2bf:last}, which requires $\mathcal{O}(nRr\log n)$ operations. Note that this computation is equivalent to the computation done for the middle-level in Algorithm $4.1$ in~\cite{liu2021butterfly} and sketching is not required for the middle level factor when the matrix is given as a low-rank decomposition. Overall, the complexity of~\cref{alg:low_rank_butterfly_rand} is $\mathcal{O}(nRr\log n + nRr^2\log n + nr^3 \log n)$. 

\begin{algorithm}
	\caption{Randomized\_LR\_to\_Butterfly}
    \label{alg:low_rank_butterfly_rand}
	\begin{algorithmic}[1]
		\State \textbf{Input:} Factors $\mat A$ and $\mat B$ of a $n\times n$ low-rank matrix with rank $R$, butterfly levels $L$ and butterly rank $r$, and oversampling factor $p$

        \State \textbf{Output:} Butterfly cores $\tsr S^{i},i=1,\ldots,L+2$.
        
            \State Reshape $\mat A$ of shape $n \times R$ to $\tsr A$ of shape $2 \times \dots \times 2 \times c \times R$ and $\mat B$ of shape $n \times R$  to $\tsr B$ of shape $2 \times \dots \times 2 \times c \times R$ by using ~\eqref{eq:indices} 
            \State{Generate random tensor $\tsr O$ of order $L+2$ and shape $2 \times \dots \times 2 \times c \times (r +p)$}

               \State $m_{i_0\dots i_L q} = \sum_{z}a_{i_0\dots i_L z} \sum_{j_0,\dots, j_L}b_{j_0\dots j_Lz}o_{j_0\dots j_L q}$\label{line:lr2bf:sketching_outer}  
                \Statex $n_{j_0\dots j_Lq} = \sum_{z}b_{j_0\dots j_L z}\sum_{i_0,\dots, i_L}a_{i_0\dots i_Lz}o_{i_0\dots i_L q}$
            \State $\mat S^{1}(i_0, \dots, i_{L-1},:,:),~~\mat V(i_0,\dots,i_{L-1}, :,:)\leftarrow$QRCP\big($\mat M(i_0, \dots, i_{L-1}, :,:)$,$r$\big)
            \Statex  $\mat S^{L+2}(j_0, \dots, j_{L-1},:,:),\mat W(j_0,\dots,j_{L-1}, :,:)\leftarrow$QRCP\big($\mat N (j_0, \dots, j_{L-1},:,:)$,$r$\big)\Comment{QR with column pivoting }\label{line:lr2bf:qrcp_outer}             
            \For{\( l = 1 \) to \( H \) }
                    \State{Generate random tensor $\tsr O$ of order $L+2-l$ and shape $2 \times \dots \times 2 \times c \times (r +p)$}
                    \State $m_{i_0 \dots i_{L-l}j_0\dots j_{l-1}r_lq} = \sum_{z}v_{i_0\dots i_{L-l}j_0 \dots j_{l-2} r_{l}z} \sum_{j_l,\dots,j_{L}}b_{j_0\dots j_Lz} o_{j_l\dots j_{L}q}$\label{line:lr2bf:sketching_inner}  
                    \Statex \hspace{1.2em} $n_{i_0\dots i_{l-1}j_0 \dots j_{L-l}r_{L+2-l}q} = \sum_{z}w_{i_0 \dots i_{l-2}j_0\dots j_{L-l} r_{L+2-l}z} \sum_{i_l,\dots,i_{L}}a_{i_0\dots i_Lz} o_{i_l\dots i_{L}q}$ 
                \State{ $\mat S^{l+1}(i_0, \dots, i_{L-l}, j_0, \dots, j_{l-1},:,:), \mat V(i_0, \dots i_{L-l}, j_0, \dots, j_{l-1})$
                \Statex \hspace{3em}$\leftarrow$ QRCP($\mat M(i_0, \dots, i_{L-l}, j_0, \dots, j_{l-1})$,$r$)}\label{line:lr2bf:qrcp_inner} 
                \Statex\hspace{1.2em} { $\mat S^{L+2-l}(i_0, \dots, i_{l-1}, j_0, \dots, j_{L-l},:,:), \mat W(i_0, \dots, i_{l-1}, j_0, \dots, j_{L-l},:,:)$
                \Statex \hspace{3em}$ \leftarrow$QRCP($\mat N( i_0, \dots, i_{l-1}, j_0, \dots, j_{L-l},:,:)$,$r$)}
            \EndFor
            \State 
                \begin{align*}
                    s^{H+1}_{i_0 \dots i_Hj_0\dots j_{H-1}r_H r_{H+1}}= &\sum_{q} s^{H+1}_{i_0 \dots i_{H-1}j_0\dots j_{H-1}r_H q}\\
                    &\sum_{z}v_{i_0\dots i_{H-1} j_0\dots j_{H-1}qz}w_{i_0 \dots i_{H-1}j_0\dots j_{H-1}r_{H+1} z} 
                \end{align*}
                \label{line:lr2bf:last}
	\end{algorithmic}
\end{algorithm}

\subsection{Full Tensor-Based Butterfly Completion}
\label{subsec:overall}
The overall workflow and cost of the proposed tensor completion algorithms for butterfly matrices are as follow. Assuming a butterfly rank $r$ and given $|\Omega|=\mathcal{O}(nr\log n)$ (which is the empirical estimate of the minimal number entries required by the butterfly completion), we first generate the initial guess  by low-rank completion with a constant number of iterations and rank estimate $R$, followed by the low-rank to butterfly conversion (\cref{alg:low_rank_butterfly_rand}). Unless otherwise stated, we use $R=r$. Hence the initial guess generation requires a total of $\mathcal O (nr^4 \log n)$ operations. Next, we run the butterfly completion (ALS-based \cref{alg:ALS_BF}, ADAM-based \cref{alg:ADAM_BF}, or ADF-based) for at most $t_{\max}$ iterations. Here $t_{\max}$ is also set to a constant. Each iteration requires $\mathcal O(|\Omega|r^4\log^2n + nr^6 \log n)=\mathcal O(nr^5\log^3n )$ operations for ALS and $O(|\Omega|r^2 \log^2n + nr^2\log n)=O(nr^3 \log^3n)$ operations for ADAM 
and ADF. As a baseline, we also implement several QTT tensor completion algorithms that are ALS,
ADAM or ALS-based. All of the above algorithms are implemented in Python and are available online \footnote{\url{https://github.com/navjo2323/butterfly-matrix-completion}}.

\textit{A remark on sampling complexity.} We remark that identification of the sampling complexity of various matrix and tensor completion algorithms is a still active research area. Ideally, the number of observed entries $|\Omega|$ should scale linearly with the number of degrees of freedom of the underlying matrix or tensor representation. While practical algorithms for completion have a long history of development, theoretical guarantees regarding their optimal sampling complexity have only recently been established. For low-rank matrix completion, earlier works can reach the optimal sampling complexity $|\Omega|=\mathcal{O}(rn)$ only for nuclear norm-based algorithms
\cite{recht2011simpler,mu2014square}, \hl{but sub-optimal complexity for ALS-based} \cite{jain2015fast} \hl{algorithms. The optimal complexity for ALS-based algorithms has only been recently proved} \cite{kelner2023matrix}. \hl{For QTT (or TT) completion, a very recent work proves the optimal $|\Omega|=\mathcal{O}(r^2\log n)$ sampling complexity for Riemannian optimization-based completion algorithms} \cite{budzinskiy2023tensor}. \hl{As our paper represents the first completion algorithms for butterfly representation, we do not attempt to establish a theoretical bound for the sampling complexity, but rather focus on empirically demonstrating the effectiveness of the proposed butterfly completion algorithms for highly oscillatory matrices. We empirically show that the sampling complexity of $|\Omega|=\mathcal{O}(rn\log n)$ is sufficient for several numerical examples in} \cref{sec:experiments} \hl{and leave its theoretical justification as a future work.}

\section{Numerical Experiments}\label{sec:experiments}
This section provides several numerical examples to demonstrate the convergence, accuracy and efficiency of the proposed completion algorithms when applied to incomplete data from high-frequency Green's function for wave equations
(\cref{sec:green}), \hl{inverse scattering applications} (\cref{sec:inverse}) \hl{and Radon transforms} (\cref{sec:Radon}). 
Specifically, the proposed algorithms are butterfly completion algorithms that are ALS-based
(\cref{alg:ALS_BF}), ADAM-based (\cref{alg:ADAM_BF}), and ADF-based (\cref{alg:butterfly-adf}). 
As for the baseline algorithms, we use the ALS-based low-rank matrix completion algorithm
(\cref{alg:ALS_BF} \hl{with $L=0$) and the QTT tensor completion algorithms that are ALS-based} (\cref{alg:ALS_QTT}), ADAM-based and ADF-based. A convergence tolerance of $\epsilon = 10^{-3}$ is used for all the algorithms. All the proposed and baseline algorithms are implemented in Python and executed using one CPU processor of the Perlmutter machine at NERSC in Lawrence Berkeley National Laboratory, which has AMD EPYC 7763 processors and 128GB of 2133MHz DDR4 memory.

\subsection{Green's function for Helmholtz equation}\label{sec:green}

In this subsection, we consider a square matrix discretized from the 3D free-space Green's function for high-frequency Helmholtz equations. Let $n = c \cdot 2^L$ denote the matrix size with $c=4$ and $L$ representing the number of butterfly or QTT levels. The matrix $\mat{\bar{T}}$ is given by
$\bar{t}_{i,j}=\frac{\mathrm{exp}(-\mathrm{i}\omega\rho)}{\rho},$
where $\omega$ is the wavenumber, $\rho=|x^{i}-y^{j}|$, $x^{{i}}=(\frac{i_1}{\sqrt{n}},\frac{i_2}{\sqrt{n}},0)$, $y^{{j}}=(\frac{j_1}{\sqrt{n}},\frac{j_2}{\sqrt{n}},1)$. Here, $i_1=\lfloor (i-1)/\sqrt{n} \rfloor$, $i_2=(i-1)\mod \sqrt{n}$, and similarly for $j_1$ and $j_2$. In this test, $\omega=\frac{\sqrt{n}\pi}{5}$, i.e., one wavelength is discretized into 10 points. We apply the completion algorithms to the permuted matrix $\mat{T} = \mat{P}\mat{\bar{T}}\mat{P}^T$ with $\mat{P}$ presenting the KD-tree-based reordering that yields the CLR property of $\mat{T}$. The completion algorithms are trained on a set of entries denoted by $\Omega$ and tested on a set of entries denoted by $\Omega_{\mathrm{test}}$, both consisting of randomly generated row and column index pairs. In other words, we define the training error using~\eqref{eq:BF_convergence}, and test error using~\eqref{eq:BF_convergence} with $\Omega\rightarrow \Omega_{\mathrm{test}}$.   

\begin{figure}[h]
	\centering
	\begin{subfigure}[b]{0.9\textwidth}
		\centering
		\includegraphics[width=\linewidth]{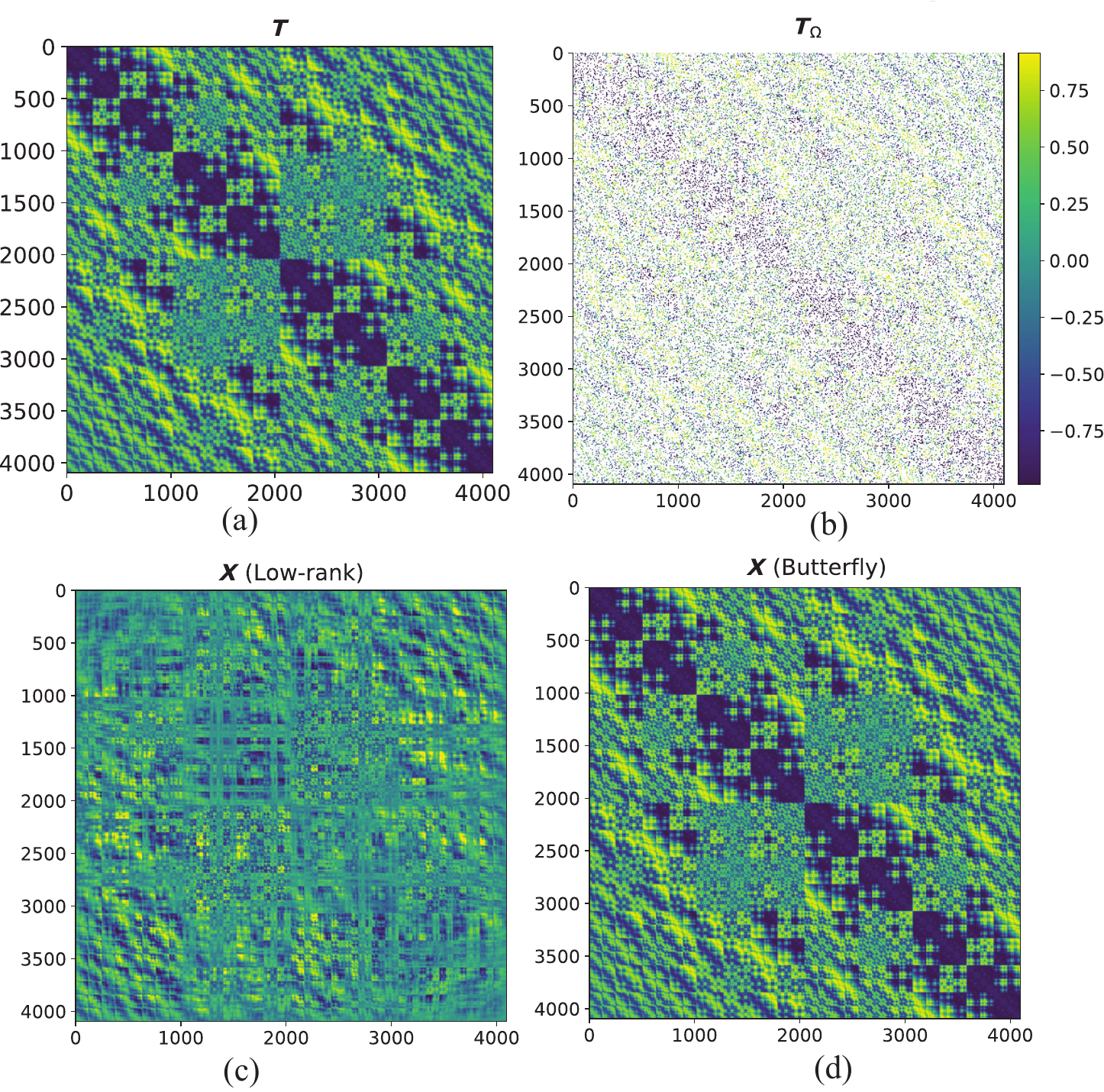}    
	\end{subfigure}
	\caption{The real part of a Green's function matrix with $n=4096$, $c=4$, $L=10$. (a) The ground truth data $\mat{T}$. (b) The observed data $\mat{T}_\Omega$ with $|\Omega|=66n\log n=3244032$. (c) The low-rank reconstruction with rank $r=400$. (d) The butterfly reconstruction using \cref{alg:ALS_BF} with rank $r=11$.\label{fig:mat_Green}}
\end{figure}

First, we consider an example of $n=4096$ and compare the reconstruction quality of the butterfly completion algorithm (Butterfly-ALS of \cref{alg:ALS_BF}) and low-rank matrix completion algorithm. The real part of the full matrix $\mat{T}$ and the observed submatrix $\mat{T}_{\Omega}$ are shown in \cref{fig:mat_Green}(a) and (b), where the number of observed entries is empirically set to $|\Omega|=66n\log n=3244032$. The low-rank completion algorithm uses rank estimate $r=400$, max iterations $t_{max}=20$, and the reconstructed full matrix $\mat{X}$ is shown in \cref{fig:mat_Green}(c). The butterfly completion algorithm (\cref{alg:ALS_BF}) uses rank estimate $r=11$, max iterations $t_{max}=20$ and use the initial guess with $R=11$. The reconstructed full matrix $\mat{X}$ is shown in \cref{fig:mat_Green}(d). As one can see, the reconstruction quality of \cref{alg:ALS_BF} is significantly better than that of the low-rank completion algorithm, given the same set of observed entries. 

\begin{figure}
	\centering
	\begin{subfigure}[b]{\textwidth}
		\centering
        {
		\includegraphics[width=\linewidth]{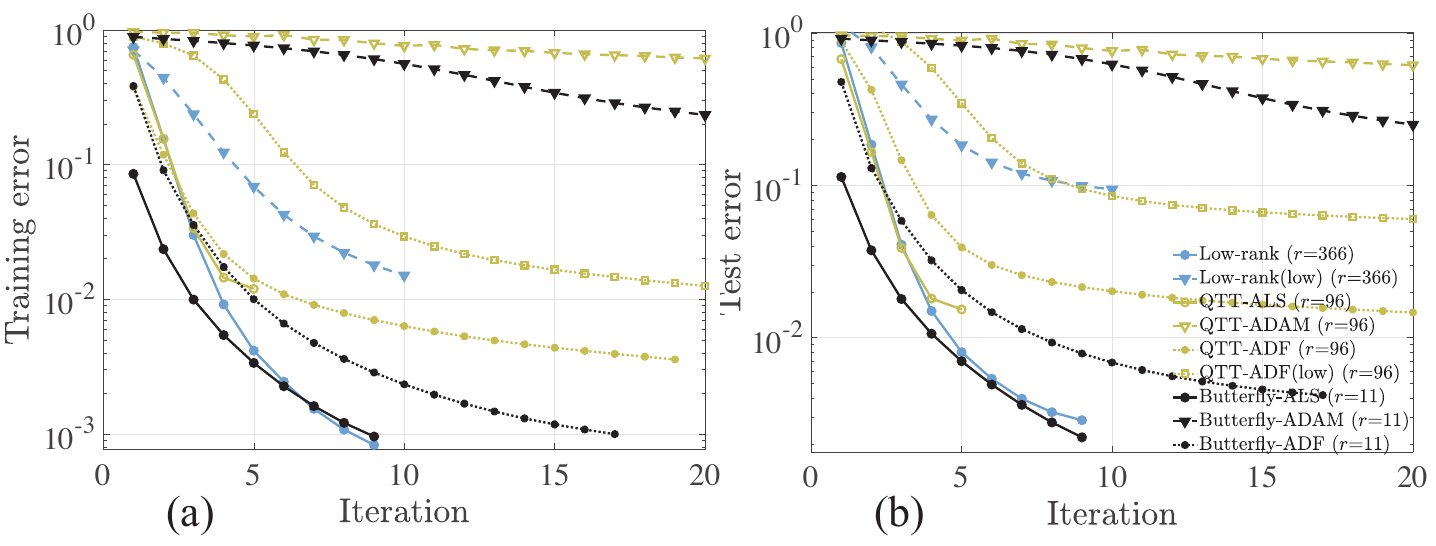}    
		}
	\end{subfigure}
	\caption{The convergence of a Green's function example with $n=16384$, $c=4$, $L=12$ using \hl{the low-rank matrix completion with $r=366$, $|\Omega|=10rn=59965440$ (i.e., $|\Omega|/n^2=0.22$), and $|\Omega|=4rn=23986176$ (i.e., $|\Omega|/n^2=0.09$, labeled as ``Low-rank(low)"), butterfly completion (ALS, ADAM, ADF) with $r=11$, $|\Omega|=6rn\log n=15138816$ (i.e., $|\Omega|/n^2=0.056$), and QTT completion (ALS, ADAM, ADF) with $r=96$, $|\Omega|=16r^2\log n=2240000$ (i.e., $|\Omega|/n^2=0.0083$). In addition, ``QTT-ADF(low)" denotes ADF-based QTT completion with $|\Omega|=8r^2\log n=1120000$. (a) Training error. (b) Test error.} \label{fig:convergence_Green_L12}}
\end{figure}

\begin{figure}
	\centering
	\begin{subfigure}[b]{\textwidth}
		\centering
    {
		\includegraphics[width=\linewidth]{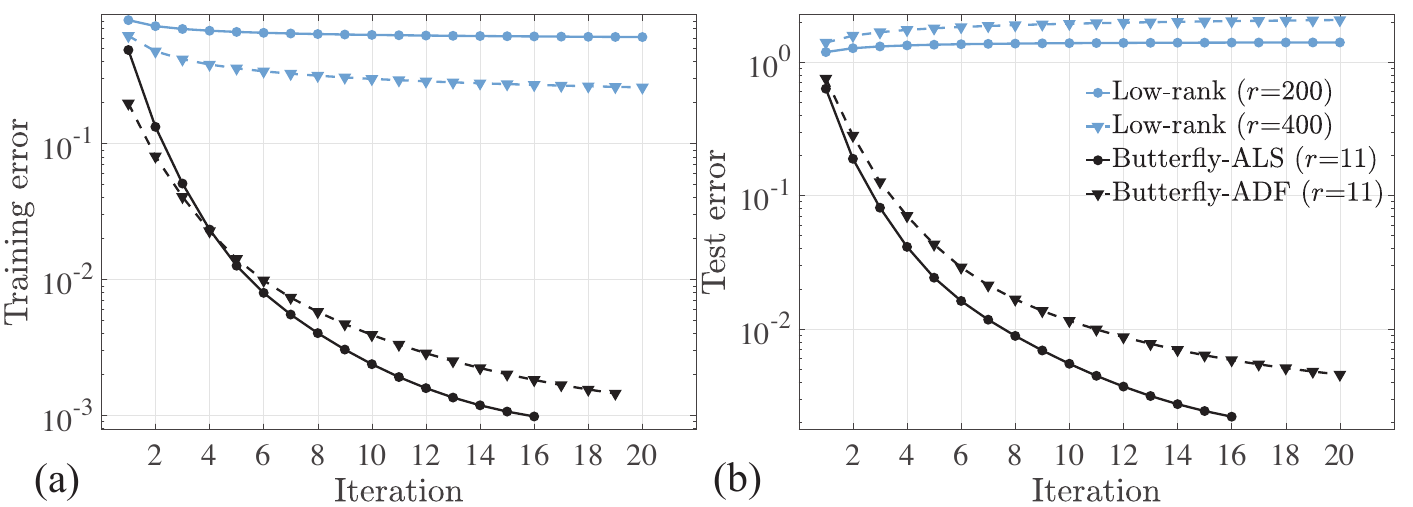}    
	}
	\end{subfigure}
	\caption{The convergence of a Green's function example with $n=65536$, $c=4$, $L=14$ using butterfly completion \hl{(ALS and ADF)} with rank $r=11$. (a) Training error with \hl{$|\Omega|=6rn\log n=69206016$}. (b) Test error with \hl{$|\Omega_{\rm test}|=rn\log n=11534336$. The results using low-rank matrix completion with rank $r=200,400$ and the same $|\Omega|$ as butterfly completion are also plotted.}\label{fig:convergence_Green_L14}}
\end{figure}

\begin{figure}
	\centering
	\begin{subfigure}[b]{\textwidth}
		\centering
                {
		\includegraphics[width=\linewidth]{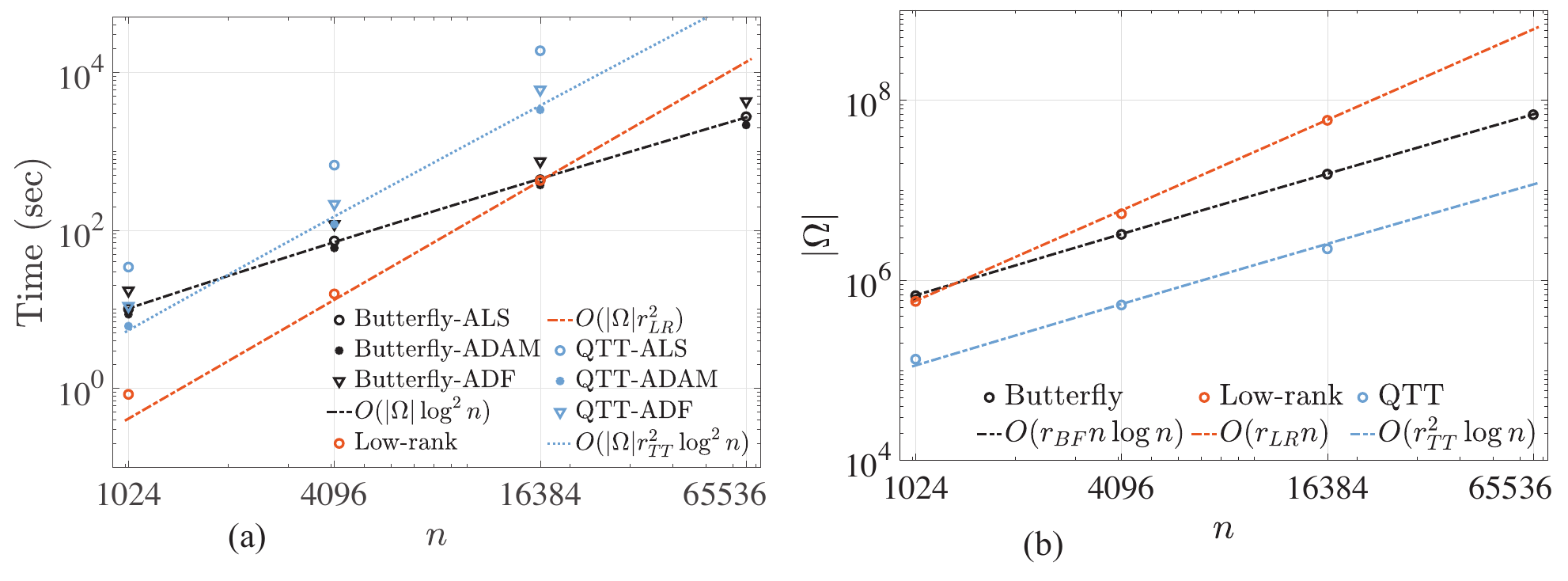}    
	}
	\end{subfigure}
	\caption{(a) Computational time per iteration and (b) $|\Omega|$ for a Green's function example with $c=4$, $L=8-14$ using the low-rank matrix completion, butterfly completion \hl{(ALS, ADAM, ADF)}, and QTT completion \hl{(ALS, ADAM, ADF)}. Butterfly completion: $|\Omega|=6rn\log n$, $r=11$ \hl{(i.e., the actual butterfly rank for $\mat T$ with tolerance $\epsilon=10^{-3}$)}. Low-rank completion: \hl{$|\Omega|=10rn$ with $r=57, 134, 366$} \hl{(i.e., the actual rank for $\mat T$ with tolerance $\epsilon=10^{-3}$)}. Note that in theory low-rank representation requires $r=\mathcal O(n)$ for sufficiently large $n$, but these numerical tests exhibit $r=\mathcal O(\sqrt{n})$ up to $n=16384$. QTT completion: \hl{$|\Omega|=16r^2\log n$ with $r=30, 54, 96$ (i.e., the actual QTT rank for $\mat T$ with tolerance $\epsilon=10^{-3}$). The QTT ranks scale approximately as $r=\mathcal O(n^{0.5})$.}\label{fig:cc_Green}}
\end{figure}

Next, we examine the convergence of the completion algorithms using the real part of the Green's function as $\mat T$. We first consider an example of $n=16384$ ($L=12$). As the ground truth, the matrix rank, butterfly rank and QTT rank with tolerance $\epsilon=10^{-3}$ are 366, 11 and 96, respectively. For all algorithms in this example, we set the rank estimates in the completion algorithms as the above true numerical ranks. Specifically, for low-rank completion, we use $r=366$ and $|\Omega|=10rn=59965440$; for the ALS-based (\cref{alg:ALS_BF}), gradient-based (\cref{alg:ADAM_BF}) and ADF-based(\cref{alg:butterfly-adf}) \hl{butterfly completion}, we use $r=11$, $R=11$ and $|\Omega|=6rn\log n=15138816$; for the ALS-based (\cref{alg:ALS_QTT}), gradient-based and ADF-based QTT completion, we set $r=96$, $R=11$ and $|\Omega|=16r^2\log n=2240000$. The convergence history for the training and test errors are shown in \cref{fig:convergence_Green_L12}. We set max iterations $t_{max}=20$ for all algorithms. 
As we set $|\Omega|$ and $r$ to approximately the minimum numbers required by each completion algorithm, all algorithms show convergence. Specifically, gradient-based methods such as Butterfly-ADAM and QTT-ADAM show slowest convergence; ADF-based algorithms show slightly worse convergence behaviors than ALS-based algorithms. For low-rank matrix completion, we also plot the convergence history for $|\Omega|=4rn=23986176$ (labeled as ``low-rank(low)"), which appears to be insufficient. For QTT-ADF, we also plot the convergence history for $|\Omega|=8r^2\log n=1120000$ (labeled as ``QTT-ADF(low)"), which also turns to be an insufficient number for convergence. Moreover, although QTT completion requires the lowest $|\Omega|$, the cost per iteration (see \cref{fig:cc_Green} \hl{for more detail) is most expensive (e.g., at $L=12$, QTT-ADF is $10\times$ slower than Butterfly-ALS, and QTT-ALS is $50\times$ slower than Butterfly-ALS) and the test errors are less satisfactory than butterfly completion.} Next, an even larger problem is tested as $n=65536$ ($L=14$). We only test the low-rank matrix completion algorithms with $r=200$ or $400$ and Butterfly-ALS/Butterfly-ADF with $r=11$. We skip Butterfly-ADAM due to its slow convergence and QTT-completion due to its high memory usage and computational cost attributed to large QTT ranks. We set $|\Omega|=69206016$ for both low-rank and butterfly completion algorithms. As can be seen from \cref{fig:convergence_Green_L14}, \hl{the low-rank matrix completion algorithms diverge due to both insufficient rank estimate and $|\Omega|$ (as will be seen in} \cref{fig:cc_Green}, it requires a projected 4 hours per iteration for low-rank completion to converge). In contrast, the butterfly algorithm converges within 20 iterations. Again, the convergence of Butterfly-ADF is slightly worse than Butterfly-ALS.
	
Next, we investigate the computational complexity of these completion algorithms with the minimum number of nonzeros it took for the algorithms to converge. As $L \in \{8, \dots,14\}$, we set the number of observed entries to \hl{$|\Omega|=10rn$} for the low-rank matrix completion algorithm, $|\Omega|=6n\log n$ for the butterfly completion algorithms, and \hl{$|\Omega|=16r^2\log n$ for the QTT algorithms.} The rank estimates are set to the empirically smallest numbers such that the completion algorithms converge for a given $n$. More specifically, we set \hl{$r=57, 134, 366$} for the low-rank algorithm, $r=11$ for the butterfly algorithms and $r=30, 54, 96$ for the QTT algorithm. Note that for $n=65536$ ($L=14$), both low-rank matrix completion and QTT completion algorithms run out of memory or are too expensive to finish. We plot the CPU time per iteration for each algorithm in \cref{fig:cc_Green}(a) and the number of observed entries $|\Omega|$ in \cref{fig:cc_Green}(b). \hl{For the time per iteration, Butterfly-ADF is about $2\times$ slower than Butterfly-ADAM as it requires two gradient-like computation per iteration (See Line} \ref{line:adf:gradient} and \ref{line:adf:reconstruct} of \cref{alg:butterfly-adf}). \hl{Surprisingly, we note that Butterfly-ALS and Butterfly-ADAM almost overlap with each other. This is likely due to the fact that the butterfly rank $r=11$ is a small constant independent of $n$. As a result, we conclude that ALS is the best algorithm among the three for butterfly completion of highly oscillatory matrices due to its fastest convergence (see} \cref{fig:convergence_Green_L12} and \cref{fig:convergence_Green_L14}) \hl{and time per iteration. On the other hand, QTT-ADF is much faster than QTT-ALS due to relatively large QTT ranks. Therefore, we conclude that ADF is the best algorithm among the three for QTT completion for highly oscillatory matrices. Note that although QTT completion requires less $|\Omega|$ than butterfly completion (see} \cref{fig:cc_Green}(b)), \hl{the time per iteration at $L=12$ is $3349$s for QTT-ADF but only $373$s for Butterfly-ALS. This time difference will be even larger for larger $L$. In terms of asymptotic behaviors, $|\Omega|$ scales as $\mathcal O(n\log n)$ for both butterfly ($r$ stays as a small constant) and QTT completion ($r=\mathcal{O}(n^{0.5})$) algorithms, and super-linearly for low-rank completion. As a result, the butterfly completion algorithms empirically achieve the lowest total asymptotic computational complexity, $\mathcal O(n\log^3n)$, among all the completion algorithms.}

\begin{figure}
	\centering
	\begin{subfigure}[b]{\textwidth}
		\centering
                {
            \includegraphics[width=\linewidth]{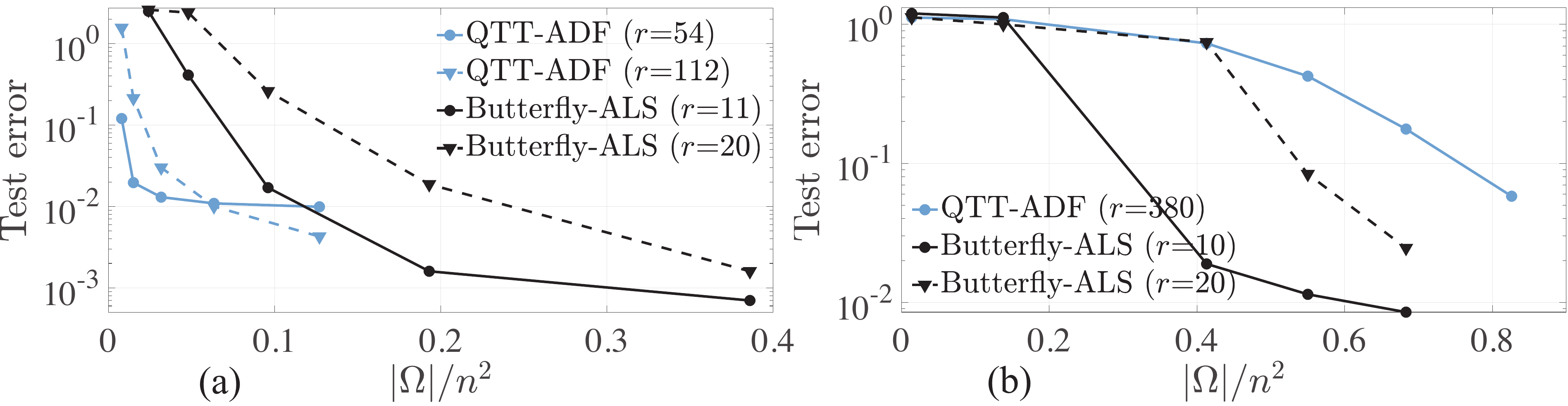}    
	}
	\end{subfigure}
	\caption{Test errors with varying observation counts $|\Omega|$ and rank estimates $r$ for butterfly completion with ALS and QTT completion with ADF. (a) A Green's function example with $c=4$, $L=10$. The actual rank for $\mat T$ with tolerance $\epsilon=10^{-3}$ is 10 for butterfly decomposition and 54 for QTT decomposition. (b) A Radon transform example with $c=4$, $L=8$. The actual rank for $\mat T$ with tolerance $\epsilon=10^{-3}$ is 10 for butterfly decomposition and 380 for QTT decomposition. We remark that for Radon transforms, QTT exhibits full ranks and requires $|\Omega|=n^2$ for accurate recovery. \label{fig:error_nnz_rank}}
\end{figure}

\hl{Moreover, we study the effect of rank estimate and $|\Omega|$ on the test error for QTT-ADF and Butterfly-ALS (i.e., the best performing algorithms out of the three QTT completion algorithms and the tree butterfly completion algorithms) using an example of $n=4096$ ($L=10$). For QTT completion, the actual rank with tolerance $\epsilon=10^{-3}$ is 54 and we set the rank estimates to $r=54,112$. For butterfly completion, the actual rank with tolerance $\epsilon=10^{-3}$ is 11 and we set the rank estimates to $r=11,20$. The test errors with varying observation ratios $|\Omega|/n^2$ are plotted in} \cref{fig:error_nnz_rank}(a). \hl{For most values of observation ratios, the use of rank estimate higher than necessary leads to over-parameterization and hence worse test errors. For this Green's function example, QTT completion typically requires lower observation ratios than butterfly completion (e.g., QTT-ADF reaches an error of 0.013 at $|\Omega|/n^2=0.03$ while Butterfly-ALS reaches an error of 0.017 at $|\Omega|/n^2=0.1$), but struggles to achieve the same level of test errors as the ratio increases.}

\subsection{Inverse wave scattering}\label{sec:inverse}
\begin{figure}
	\centering
	\begin{subfigure}[b]{\textwidth}
		\centering
		\includegraphics[width=\linewidth]{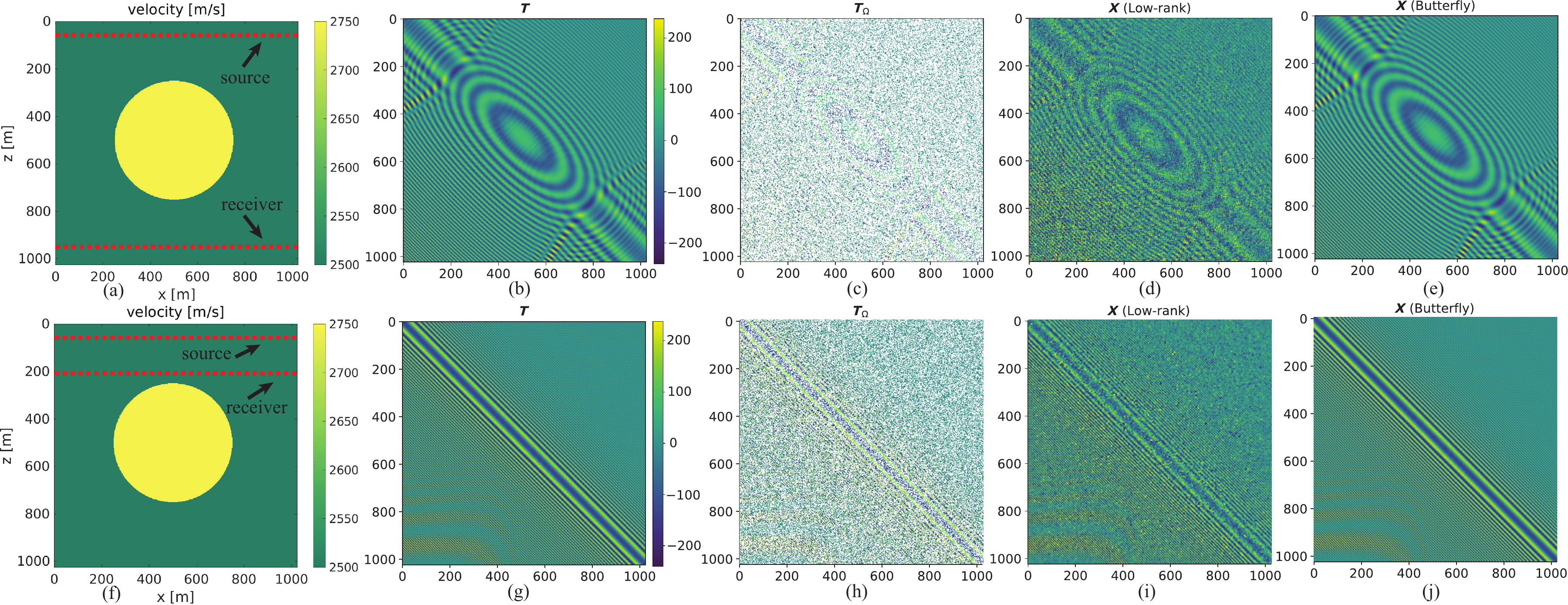}    
	\end{subfigure}
	\caption{Top: the transmission problem. Bottom: the reflection problem. From left to right: the velocity model under consideration, the full target matrix $\mat{T}$, the observed matrix entries, the reconstructed matrix $\mat{X}$ by the low-rank completion algorithm and the reconstructed matrix $\mat{X}$ by the butterfly completion algorithm. Here we set $n=1024$, $c=4$, $L=8$, $|\Omega|=27n\log n$ for the transmission problem and $|\Omega|=48n\log n$ for the reflection problem. \label{fig:mat_seismic}}
\end{figure}

Next, we consider a more realistic inverse scattering application governed by the following 2D Helmholtz equation
\begin{equation}
\nabla^2 p({\mathbf x})+\frac{\omega^2}{v^2({\mathbf x})}p({\mathbf x})=-f({\mathbf x}),\label{eq:visco-acoustic}
\end{equation}
in a physical domain of $1000 \mathrm{~m} \times 1000\mathrm{~m}$. Here ${\mathbf x}=(x_1,x_2)$, $f({\mathbf x})$ is the excitation function, $p(\mathbf x)$ is the pressure wave field, $\omega$ is the angular frequency, and $v({\mathbf x})$ is the velocity profile. Here $v({\mathbf x})=2500 \mathrm{~m/s}$ everywhere except for a circular target region where $v({\mathbf x})=2750 \mathrm{~m/s}$ (see the first column of \cref{fig:mat_seismic}). In the inverse scattering application, one is interested in generating an array of point sources $f({\mathbf x})=\delta(\mathbf x-\mathbf{x}^s_i)$ at locations $\mathbf{x}^s_i$, $i=1,2,\ldots, n$ and measuring the pressure wave fields $p(\mathbf{x}^o_j)$ at locations $\mathbf{x}^o_j$, $j=1,2,\ldots, n$; from this measured $n\times n$ Green's function matrix, one can infer the unknown velocity ($v({\mathbf x})=2750 \mathrm{~m/s}$ is the ground truth to be inferred) inside the target region. Here, we demonstrate how to reconstruct the $n\times n$ matrix in some compressed formats from incomplete measurement, which can be potentially used to infer the unknown velocity profile. We consider both the transmission problem where the source and observer arrays (with equally spaced points) are located at $z=50$ m and $z=950$ m (see the top panel of \cref{fig:mat_seismic}), and the reflection problem where the source and observer arrays are located at $z=50$ m and $z=200$ m (see the bottom panel of \cref{fig:mat_seismic}). 

To generate the Green's function matrix $\mat{T}$, we leverage the simulation package WAVEFORM \cite{WAVEFORM} to simulate the responses $p(\mathbf{x}^o_j)$, $j=1,2,\ldots, n$ for each excitation function $f({\mathbf x})=\delta(\mathbf x-\mathbf{x}^s_i)$. The frequency is set to $\omega=500\pi$ and $n$ is set to 1024 (i.e., $c=4$, $L=8$). Once the matrix is fully generated, we assume only $|\Omega|$ entries of the matrix are observed and our goal is to fully recover $\mat{T}$. For the transmission problem, we empirically set $|\Omega|=27n\log n$, $r=9$ for the butterfly completion \cref{alg:ALS_BF} and $r=120$ for the low-rank completion. For the reflection problem, we empirically set $|\Omega|=48n\log n$, $r=12$ for the butterfly completion \cref{alg:ALS_BF} and $r=120$ for the low-rank completion. The reconstructed matrices are plotted in the last two columns of \cref{fig:mat_seismic}. Clearly, the butterfly completion algorithm achieves significantly better reconstruction quality than the low-rank completion algorithm.

\subsection{Radon transform}\label{sec:Radon}

In this subsection, we consider a matrix $\mat{T}$ discretized from a 1D Radon transform:
\begin{eqnarray*}
	t_{i,j}=\mathrm{exp}(2\pi\mathrm{i}\phi(x^{{i}},y^{{j}}))\\
	\phi(x,y)=x\cdot y +  c|y|  \label{eq:randon2d}\\
	c=(2+\sin(2\pi x))/8
\end{eqnarray*}
with $x^{{i}} = \frac{i}{n}$, $y^{{j}} = j-\frac{n}{2}$, $i,j=1,2,\ldots,n$. This matrix is known to exhibit full rank when represented as low-rank or QTT algorithms, but show constant rank when represented as butterfly.  

\begin{figure}
	\centering
	\begin{subfigure}[b]{0.9\textwidth}
		\centering
		\includegraphics[width=\linewidth]{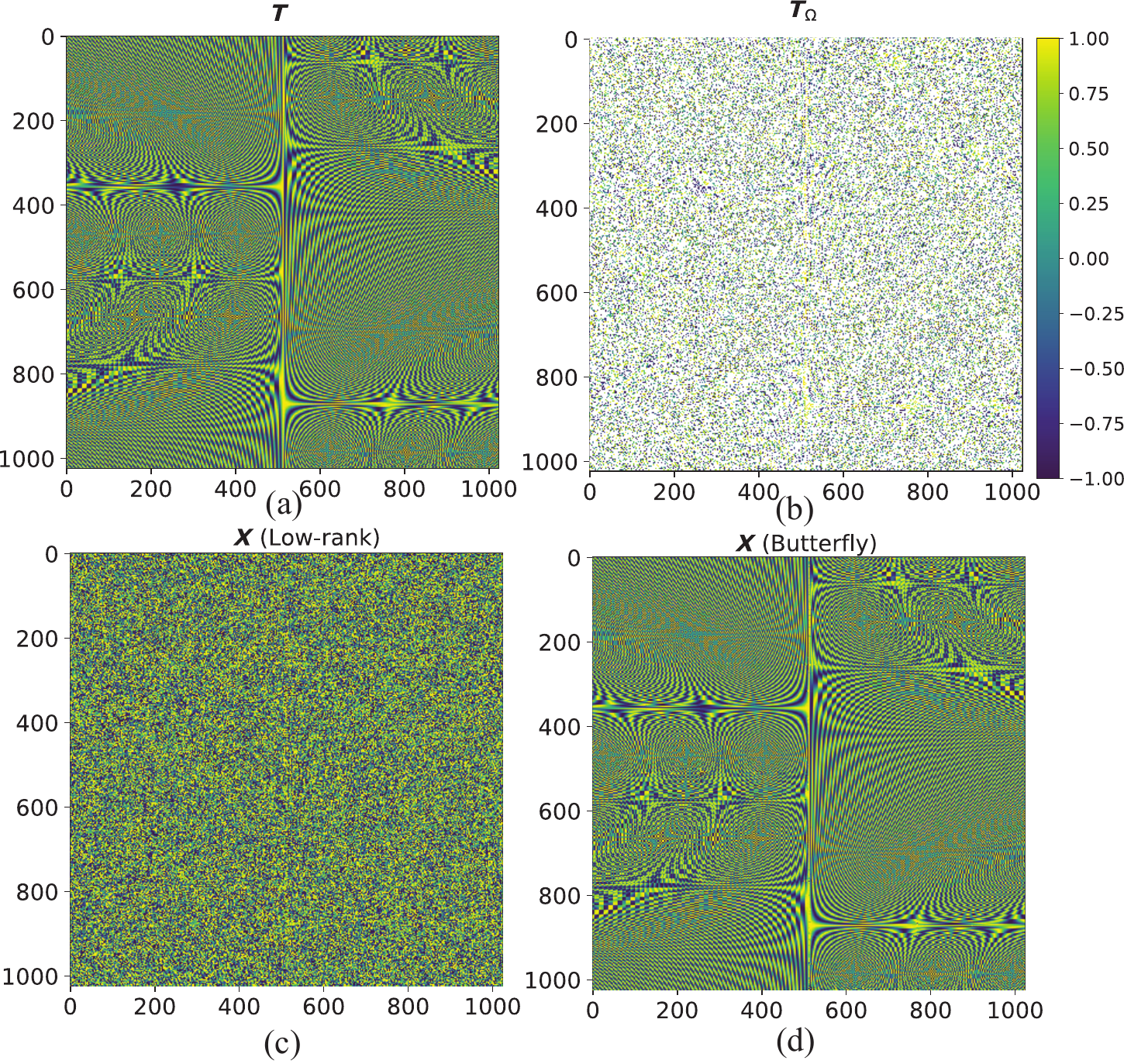}    
	\end{subfigure}
	\caption{The real part of a Radon transform matrix with $n=1024$, $c=4$, $L=8$. (a) The ground truth data $\mat{T}$. (b) The observed data $\mat{T}_\Omega$ with $|\Omega|=30n\log n=307200$. (c) The low-rank reconstruction using low-rank completion with rank $r=400$. (d) The butterfly reconstruction using \cref{alg:ALS_BF} with rank $r=10$. Note that only the real part of the data is plotted. \label{fig:mat_Radon}}
\end{figure}

First, we consider an example of $n=1024$ and compare the reconstruction quality of the butterfly completion algorithm (\cref{alg:ALS_BF}) and low-rank matrix completion algorithm. The real part of the full matrix $\mat{T}$ and the observed submatrix $\mat{T}_{\Omega}$ are shown in \cref{fig:mat_Radon}(a) and (b), where the number of observed entries is empirically set to $|\Omega|=30n\log n=307200$. The low-rank matrix completion algorithm uses rank estimate $r=400$, max iterations $t_{max}=20$, and the reconstructed full matrix $\mat{X}$ is shown in \cref{fig:mat_Radon}(c). The butterfly completion algorithm (\cref{alg:ALS_BF}) uses rank estimate $r=10$, max iterations $t_{max}=20$ and the initial guess with $R=10$. The reconstructed full matrix $\mat{X}$ is shown in \cref{fig:mat_Radon}(d). As one can see, the reconstruction quality of \cref{alg:ALS_BF} is significantly better than that of the low-rank completion algorithm, given the same set of observed entries. This is no surprise as the low-rank completion algorithm requires $r= \mathcal O (n)$ and $|\Omega|= \mathcal O (n^2)$ for a good reconstruction. 

\begin{figure}
	\centering
	\begin{subfigure}[b]{\textwidth}
		\centering
                {
		\includegraphics[width=\linewidth]{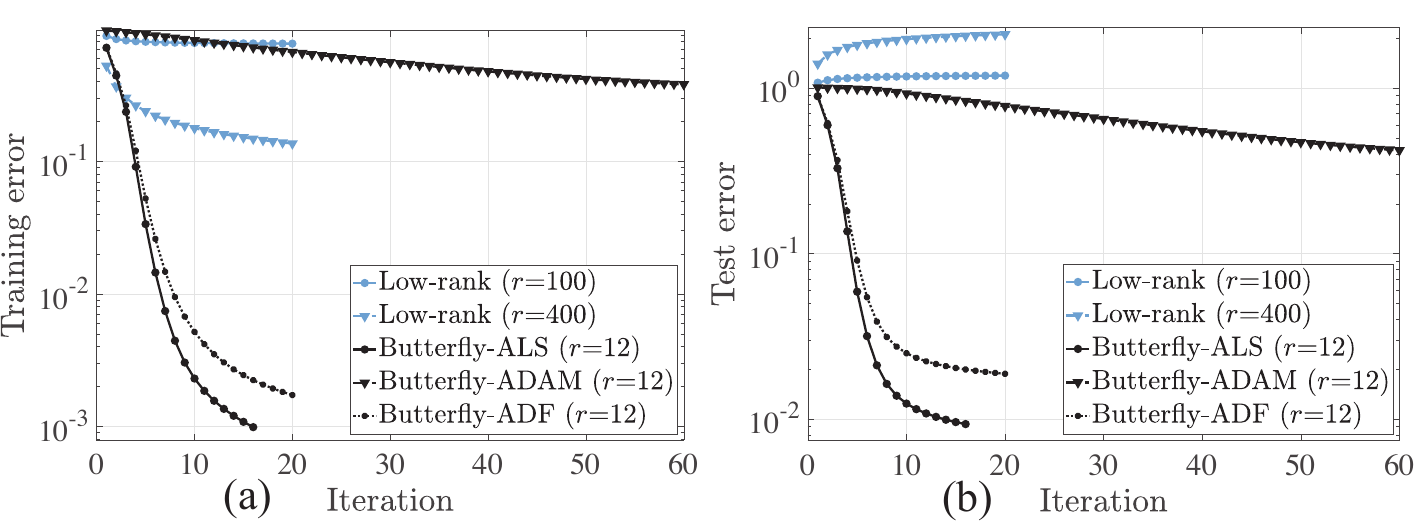}    
	}
	\end{subfigure}
	\caption{The convergence of a Radon transform example with $n=4096$, $c=4$, $L=10$ using butterfly completion \hl{(ALS, ADAM, ADF)} with rank $r=12$ \hl{(the actual butterfly rank for $\mat T$ with tolerance $\epsilon=10^{-3}$ is 10)}, \hl{$|\Omega|=7rn\log n=3440640$ and $|\Omega_{\rm test}|=rn\log n=491520$. (a) Training error. (b) Test error.} \hl{The results using low-rank matrix completion with rank $r=100,400$ and the same $|\Omega|$ as butterfly completion are also plotted. Note that the actual rank for $\mat T$ with tolerance $\epsilon=10^{-3}$ is $r=4096$ (i.e., full rank) for low-rank decomposition.}\label{fig:convergence_Radon}}
\end{figure}

\begin{figure}
	\centering
	\begin{subfigure}[b]{\textwidth}
		\centering
                {
		\includegraphics[width=\linewidth]{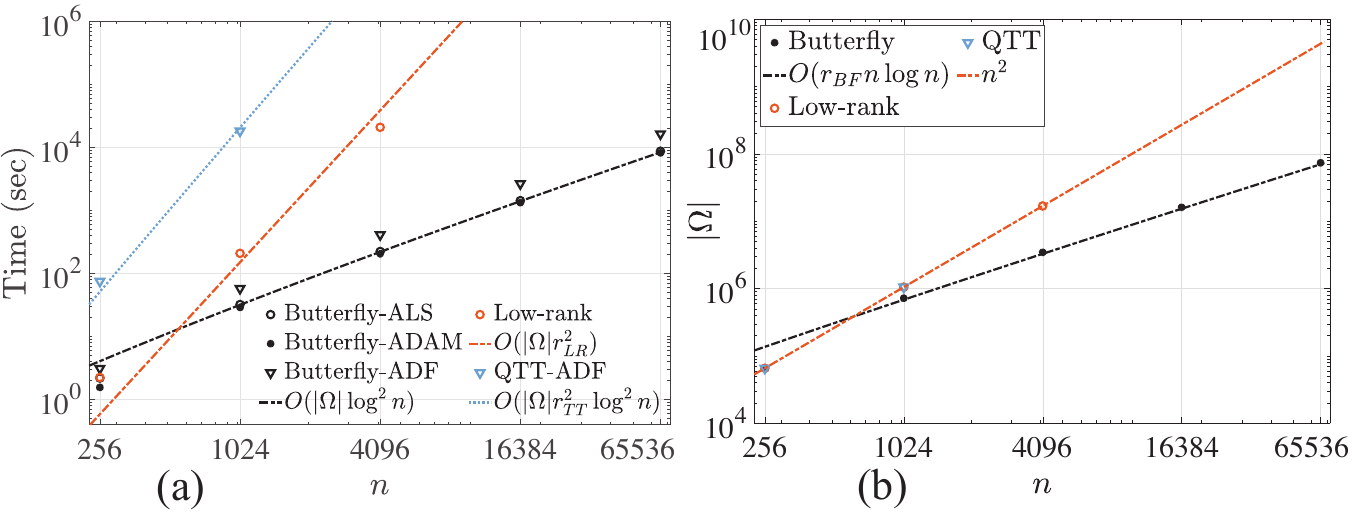}  
	}  
	\end{subfigure}
	\caption{(a) Computational time per iteration and (b) $|\Omega|$ for the Radon transform example with $c=4$, $L=6-14$ using the low-rank matrix completion, butterfly completion \hl{(ALS, ADAM, ADF) and QTT completion (ADF).} Butterfly completion: $|\Omega|=6n\log n$, $r=10$. Low-rank completion: $|\Omega|=n^2$ , $r=n$ (this requires all the entries as $\mat T$ is full-rank). \hl{QTT completion: $|\Omega|=n^2$ (this requires all the entries as the actual TT-rank for $\mat T$ is $r=140, 380, 1140$ for $L=6,8,10$).}\label{fig:cc_Radon}}
\end{figure}

Next, we examine the convergence of the completion algorithms. We first consider an example of $n=4096$ ($L=10$) using the ALS-based (\cref{alg:ALS_BF}) and gradient-based (\cref{alg:ADAM_BF}) \hl{and ADF-based} (\cref{alg:butterfly-adf}) butterfly completion with rank $r=12$ \hl{(the actual butterfly rank is 10 with tolerance $\epsilon=10^{-3}$)}. Note that the QTT completion algorithm is not considered as the full QTT rank will yield the algorithm prohibitively expensive (\hl{see later discussion of} \cref{fig:cc_Radon} and \cref{fig:error_nnz_rank}(b) \hl{for more detail}). The number of observed entries is empirically set to $|\Omega|=7rn\log n$ for training and $|\Omega_{\mathrm{test}}|=n\log n$ for testing. \hl{As the reference, we also tested the low-rank matrix completion with rank $r=100$ or $200$, but using the same $|\Omega|$ as butterfly completion.} The convergence history for the training and test errors are shown in \cref{fig:convergence_Radon}. We set max iterations $t_{max}=20$ for all the ALS-based \hl{and ADF-based} algorithms and $t_{max}=60$ for the gradient-based algorithm. As expected, the low-rank matrix completion algorithms do not converge \hl{due to both insufficient numbers of $r$ and $|\Omega|$}. In contrast, the butterfly completion algorithms converge with a small rank estimate $r=12$. Specially, the ALS-based butterfly algorithm converges in 16 iterations, \hl{and the ADF-based butterfly algorithm is slightly worse}, while the gradient-based butterfly algorithm takes more than 60 iterations to converge. 

Next, we investigate the computational complexity of \hl{the low-rank matrix, QTT and butterfly} completion algorithms for the radon transform matrices. As $L\in \{6,\dots,14\}$, we set the number of observed entries to $|\Omega|=n^2$ \hl{for low-rank completion and QTT completion} as they need all the entries to be observed to converge to sufficient accuracy, and $|\Omega|=6rn\log n$ for the butterfly algorithms. \hl{For QTT completion, we only test QTT-ADF as the actual QTT ranks with tolerance $\epsilon=10^{-3}$) for $L=6,8,10$ are $140, 380, 1140$, which are almost full-rank and scale as $\mathcal{O}(n)$. The rank estimate is $r=10$ for the butterfly algorithms, $r=n$ for the low-rank algorithm, and $r=140,380$ for the QTT completion algorithm.} We plot the CPU time per iteration for each algorithm in \cref{fig:cc_Radon}(a) and the number of observed entries $|\Omega|$ in \cref{fig:cc_Radon}(b). \hl{Similar to the Green's function example, the time per iteration for ALS-based and ADAM-based butterfly algorithms almost overlap with each other, and ADF-based butterfly algorithm is about $2\times$ slower. Clearly, the butterfly algorithms achieve the quasi-linear asymptotic complexity as $\mathcal O(n\log^3n)$. In stark contrast, low-rank and QTT completion algorithms scale at least as $\mathcal{O}(n^4)$ due to the unfavorable $|\Omega|=n^2$ and $r=\mathcal{O}(n)$ scaling. Note that at $n=1024$, it requires 5 hours per iteration for QTT completion. Obviously, if all the entries of the matrix are available, one should use low-rank or QTT decomposition instead of low-rank or QTT completion.}

\hl{Finally, we study the effect of rank estimate and $|\Omega|$ on the test error for QTT-ADF and Butterfly-ALS using an example of $n=1024$ ($L=8$). For QTT-ADF, the actual rank with tolerance $\epsilon=10^{-3}$ is 380 and we set the rank estimates to $r=380$. For Butterfly-ALS, the actual rank with tolerance $\epsilon=10^{-3}$ is 10 and we set the rank estimates to $r=10,20$. The test errors with varying observation ratios $|\Omega|/n^2$ are plotted in} \cref{fig:error_nnz_rank}(b). \hl{For most values of observation ratios, the use of rank estimate higher than necessary leads to over-parameterization and hence worse test errors. It is clear that QTT-ADF requires almost all matrix entries and hence does not represent an effective completion algorithm for Radon transform matrices.}

\section{Conclusion}

In this work, we proposed \hl{three} low-complexity butterfly completion algorithms, respectively based on ALS, ADAM \hl{and ADF}, for reconstructing highly oscillatory operators discretized as $n\times n$ matrices requiring only $\mathcal{O}(n\log n)$ observed matrix entries and $\mathcal{O}(n \log^3 n)$ computational cost empirically. The proposed algorithms are based on the reformulation of butterfly matrices into higher order tensors, which allows for efficient implementation via tensor computations. Numerical results demonstrate the superior accuracy, convergence, and efficiency of the proposed algorithms for highly oscillatory operators such as Radon transforms and Green's functions in seismic applications, when compared with other completion algorithms such as low-rank matrix completion and QTT completion. \hl{While our experiments utilized only synthetic datasets, the oscillatory operators studied here, such as high-frequency Green's functions and Radon transforms, form the mathematical foundation of several critical imaging tasks, where physical or cost constraints limit data acquisition. The proposed butterfly completion algorithm is potentially useful for addressing missing-trace problem for high-resolution reverse time migration of seismic exploration, the frequency-domain gaps inherent in low-dose medical CT, and the ``missing wedge" problem in electron tomography.} 

Future work includes exploiting butterfly completion for stochastic low-rank decompositions, high-performance and parallel implementations of the butterfly contraction operations, as well as \hl{integrating the tensorized reformulation of the butterfly decomposition} into neural network architectures to facilitate efficient learning with butterfly matrices.

\section*{Acknowledgements}

This research has been supported by the U.S. Department of Energy, Office of Science, Office of Advanced Scientific Computing Research, Mathematical Multifaceted Integrated Capability
Centers (MMICCs) program, under Contract No. DE-AC02-05CH11231 at Lawrence Berkeley National Laboratory. This research used resources of the National Energy Research Scientific Computing Center, a DOE Office of Science User Facility supported by the Office of Science of the U.S. Department of Energy under Contract No. DE-AC02-05CH11231 using NERSC award ASCR-ERCAP0024170.

\appendix
\section{Efficient Implementation of ALS for Butterfly Completion}
\label{sec:appendix}
Implementation of the ALS algorithm for butterfly completion as described in~\cref{alg:ALS_BF} efficiently can be a challenging task when the number of levels of the butterfly decomposition increase. This is due to computations with high-order parameter tensors with their order depending on the number of levels in the butterfly decomposition. In this section, we describe a method to efficiently implement ALS for butterfly completion and QTT completion with a single software kernel that involves sparse tensor computations. This method also allows for efficient sampling of entries of a butterfly matrix (Section 3.4 of \cite{liu2021sparse}). \hl{Our implementation is in Python and is open-source~\footnote{\url{https://github.com/navjo2323/butterfly-matrix-completion}}}.

\subsection{ALS for QTT completion}
\label{subsec:ALS_QTT}
Firstly, we present ALS for QTT format in~\cref{alg:ALS_QTT}. Similar to butterfly completion,~\cref{alg:ALS_QTT} proceeds by taking the observed input tensor, rank $r$ (assumed to be the same for all tensors here, but can be easily generalized to different ranks), initial guess for the tensors, and a parameter for controlling the maximum number of iterations if the convergence criterion is not satisfied. For the outermost tensors, $\tsr S^1$ and $\tsr S^{L+1}$, the algorithm solves for each fiber (potentially in parallel) requiring an $r \times r$ positive semi-definite (PSD) solve. The steps to compute the left and right-hand sides of the systems are similar to \cref{alg:ALS_BF}, however, the $L$ matrix-vector products are indexed differently and are given in line~\ref{line:QTT:compute_u} and line~\ref{line:QTT:compute_v}. Updating the outermost tensors requires at most $L$ matrix-vector products for each nonzero with a total cost $\mathcal{O}(|\Omega|Lr^2)$ (including the outer product) 
and solving for each fiber of the tensor requires $O(r^3)$ as there are only $4$ fibers to be solved for.
For the inner tensors, $\big \{ \tsr S^i : i \in \{2, \dots, L\} \big\}$, the algorithm solves for each slice requiring $r^2 \times r^2$ PSD solve. The left and right-hand sides of the systems are computed in a similar manner as in \cref{alg:ALS_BF} except for the matrix-vector products which are given in line~\ref{line:QTT:compute_u_inner} and line~\ref{line:QTT:compute_v_inner}. Similar to \cref{alg:ALS_BF}, the result is reshaped into a matrix in row-major ordering to update the inner tensors. The total cost of performing the computations in line~\ref{line:QTT:compute_u_inner} and line~\ref{line:QTT:compute_v_inner} is $\mathcal O(|\Omega|Lr^2)$, and the total cost of forming the left-hand sides $\mathcal O \big (|\Omega| (Lr^2 +r^4) \big)$, which includes the outer products of Kronecker products in line~\ref{line:QTT:Kronecker_outer}. The solve can be performed in a total cost of $\mathcal O(r^6)$. Since the number of tensors to be solved is $L$, the total cost of one iteration of~\cref{alg:ALS_QTT} is $\mathcal O \big(|\Omega| L (Lr^2 + r^4 ) + Lr^6 \big)$.

It should be noted that the~\cref{alg:ALS_QTT} is exactly the same as~\cref{alg:ALS_BF} except for the indexing of the matrix-vector products used to form the left and right-hand sides. Therefore, by introducing an indexing on the input tensor entries such that the outer tensors are reshaped into matrices and the inner tensors are reshaped into order $3$ tensors, one can use exactly the same software kernel for tensor train, QTT and butterfly completion. This also means that the computational benefits of storing the sparse tensor entries in the  compressed sparse fiber (CSF) format~\cite{smith2015tensor} to amortize the matrix-vector products involved in each solve can be translated to the butterfly completion algorithms. We describe this indexing in the next subsection.

\begin{algorithm}[H]
	\caption{\small \textsf{ALS for QTT Completion}\label{alg:ALS_QTT}}
    \small
	\begin{algorithmic}[1]
		\State \textbf{Input:} Observed matrix $\tsr{T}_\Omega$, rank \( r \), initial guess $\tsr S^{i},i=1,\ldots,L+1$, max iterations \( t_{\max} \).
		\State \textbf{Output:} Updated QTT cores $\tsr S^{i},i=1,\ldots,L+1$.
		\For{\( t = 1 \) to \( t_{\max} \)}
		\For {$l = 1$ to $L + 1$} 		
		\If {$l=0$} 
		\For {each $(i_0, j_{0})$}  
		\State Initialize $\vcr y = \vcr 0$, $\mat K = \mat 0$
		\For {each $(i_1, \dots, i_{L},j_1, \dots, j_{L}) \in \Omega_{i_0, j_{0}}$}
		\State $\vcr{v} = 
		\Bigg(\prod_{m=1}^{L-1} 
		\mat{S}^{m+1}(i_{m}, j_{m},:,:) 
		\Bigg)  
		\vcr{s}^{L+1}(i_L, j_L,:)$ \label{line:QTT:compute_v}
		\State Update $\vcr y \mathrel{+}= t_{i_0 \dots i_L j_0 \dots j_L} \vcr{v}$ 
		\State Update $\mat K \mathrel{+}= \vcr{v}\vcr{v}^T$
		\EndFor
		\State Solve $\mat K \vcr s = \vcr y$ and update $\vcr s^{1}(i_0, j_{0},:) \leftarrow \vcr s$    \Comment{$\mat{K} \in \mathbb{R}^{r \times r} $}
		\EndFor				
		\ElsIf{$l=L+1$} 
		\For {each $(i_L, j_{L})$}  
		\State Initialize $\vcr y = \vcr 0$, $\mat K = \mat 0$
		\For {each $(i_0, \dots, i_{L-1},j_0, \dots, j_{L-1}) \in \Omega_{i_L, j_{L}}$}
		\State $\vcr{u} = 
		\vcr{s}^{1}(i_L, j_L,:){}^T\Bigg(\prod_{m=1}^{L-1} 
		\mat{S}^{m+1}(i_{m}, j_{m},:,:) 
		\Bigg)$ \label{line:QTT:compute_u}
		\State Update $\vcr y \mathrel{+}= t_{i_0 \dots i_L j_0 \dots j_L} \vcr{u}$ 
		\State Update $\mat K \mathrel{+}= \vcr{u}\vcr{u}^T$
		\EndFor
		\State Solve $\mat K \vcr s = \vcr y$ and update $\vcr s^{L+1}(i_L, j_{L},:) \leftarrow \vcr s$ \Comment{$\mat{K} \in \mathbb{R}^{r \times r} $}
		\EndFor				
		\Else
		\For {each $(i_{l-1}, j_{l-1})$} 
		\State Initialize $\vcr y = \vcr 0$, $\mat K = \mat 0$
		\For {each $(i_{0}, \dots, i_{l-2},i_{l},\dots,i_{L},j_{0}, \dots, j_{l-2},j_{l},\dots,j_{L})$}
		\State $\vcr{u} = 
		\vcr{s}^{1}(i_0, j_0,:){}^T\Bigg(\prod_{m=1}^{l-2} 
		\mat{S}^{m+1}(i_{m}, j_{m},:,:)
		\Bigg)$  \label{line:QTT:compute_u_inner}
		\State $\vcr{v} = 
		\Bigg(\prod_{m=l}^{L-1} 
		\mat{S}^{m+1}(i_{m}, j_{m},:,:) 
		\Bigg)  
		\vcr{s}^{L+1}(i_L, j_L,:)$ 
         \label{line:QTT:compute_v_inner}
		\State Update $\vcr y \mathrel{+}= t_{i_0 \dots i_L j_0 \dots j_L} (\vcr{v}\otimes\vcr{u})$ 
		\State Update $\mat K \mathrel{+}= (\vcr{v}\otimes\vcr{u})(\vcr{v}\otimes\vcr{u})^T$
        \label{line:QTT:Kronecker_outer}
		\EndFor
		\State Solve $\mat K \vcr s= \vcr y$   \Comment{$\mat{K} \in \mathbb{R}^{r^2 \times r^2} $}
		\State Update $\mat S^{l}(i_{l-1}, j_{l-1},:,:) := \mathrm{reshape}(\vcr s,r,r)$
		\EndFor		
		\EndIf
		\EndFor	
		\EndFor	
	\end{algorithmic}
\end{algorithm}

\subsection{Indexing for butterfly completion algorithms}
\label{subsec:Indexing_in_BF}
To reshape the higher-order tensors in tensorized butterfly factorization such that the outermost tensors are represented as matrices and the inner tensors are represented as order $3$ tensors, we use the inverse of the function defined in~\eqref{eq:function_ind}. Precisely, we use a function $\psi^{-1}(\cdot)$ that takes a $p$-bit binary number as input and outputs an integer given as

\begin{align}
    \psi^{-1}(b_0,\dots,b_{p-1}) =\sum_{i=0}^{p-1} b_i \, 2^{p-1-i}.
\end{align}

Therefore, a new list of factors $\mathcal Z = \big \{ \tsr Z^l: l \in \{1,\dots, L+2\} \big\}$
can be formulated based on the above reshaping, by reshaping the $\tsr S^l$ introduced in ~\eqref{eq:tensor_form_factors} as
\begin{align}
    \vcr z^1 \big( k_1,: \big) &= \vcr s^1 (i_0,\dots,i_L,:) \nonumber\\
    \vcr z^{L+2} \big( k_{L+2},: \big) &=\vcr{s}^{L+2}(j_0, \dots, j_L,:) \nonumber \\
    \mat Z^{m+1} \big( k_{m+1}, :,: \big) &= \mat S^{m+1}(i_0, \dots, i_{L-m} , j_0, \dots, j_{m-1},:,:),
\label{eq:reshaped_tensors}
\end{align}
where $m \in \{1,\dots,L\}$, $k_1=\psi^{-1}(i_0,\dots, i_L)$, $k_{L+2}=\psi^{-1}(j_0,\dots, j_L)$, and $k_{m+1} = \psi^{-1}(i_0, \dots i_{L-m} , j_0, \dots, j_{m-1})$.

By using order 2 and order 3 tensor representations $\mathcal Z$ that adhere to the tensor train format, any tensor-train completion algorithm can be used to optimize the factors. Note that similar reshapes can be defined for the QTT format as well to utilize the implementation of tensor-train completion algorithms. We use the above reshapes for implementing~\cref{alg:ALS_BF} and~\cref{alg:ALS_QTT} with a tensor-train completion optimization algorithm. In our implementation of ALS for tensor-train completion, we transpose the sparse input tensor according to the index being solved for. This results in increasing the computational efficiency by using BLAS-$3$ instead of BLAS-$2$ operations to form the left-hand sides~(see~\cite{singh2022distributed}) without constructing the CSF trees. The software kernel can further be optimized by using the computational trees which would allow for further increase in computational efficiency by reusing the matrix-vector products needed to form left and right-hand sides (see~\cite{smith2016exploration}).

\subsection[Orthogonalization sweeps for butterfly ADF]{\hl{Orthogonalization sweeps for butterfly ADF}}
\label{app:butterfly-orth}

To apply ADF~\cite{grasedyck2013alternating, grasedyck2015variants} to the 
butterfly factors $\mathcal{S}$, we require orthogonalization sweeps that 
produce orthonormal matrices in the local least-squares systems. In 
the higher-order tensor representation 
$\tsr{S}^{l}$~\eqref{eq:tensor_form_factors}, orthogonalization at level 
$l+1$ proceeds by folding the outermost mode that is not shared with the 
adjacent factor into the rank dimension, forming a tall (or wide) matrix, 
and computing its QR (or LQ) factorization. The resulting triangular factor 
$\mat{R}$ (or $\mat{\mathcal{L}}$) is then absorbed into the adjacent factor 
so that the overall product is preserved. In the reshaped representation 
where we only have order $2$ and $3$ 
tensors~\eqref{eq:reshaped_tensors}, the modes of $\tsr{Z}^{m+1}$ are 
encoded in the binary digits of the index 
$k_{m+1} = \psi^{-1}(i_0, \dots, i_{L-m}, j_0, \dots, j_{m-1})$, so 
folding a mode corresponds to pairing matrices that differ at a specific 
bit position in $k_{m+1}$ and concatenating them along the appropriate rank 
axis. We describe the QR factorization and absorption steps below for a 
left-to-right sweep; the right-to-left sweep is analogous and can be 
implemented by transposing each factor, performing the same QR-based 
procedure in reverse order, and transposing back.

\paragraph{QR factorization}
For the outer factor $\mat{Z}^1$, we perform independent reduced QR 
factorizations of each of the $2^L$ blocks of size $(c \times r)$, 
yielding orthonormal blocks $\mat{Q}_{k_1}$ and upper-triangular remainders 
$\mat{R}_{k_1}$. Here recall that $k_1=\psi^{-1}(i_0,\dots, i_L)$. For an inner factor $\tsr{Z}^{m+1}$ with 
$m \in \{1, \dots, L\}$, the mode to be folded is $i_{L-m}$, since this is 
the outermost row index that appears in $\tsr{Z}^{m+1}$ but not in the 
next factor $\tsr{Z}^{m+2}$ (compare their index sets 
in~\eqref{eq:reshaped_tensors}). In the encoding 
$k_{m+1} = \psi^{-1}(i_0, \dots, i_{L-m}, j_0, \dots, j_{m-1})$, the 
index $i_{L-m}$ occupies bit position $m$ from the least significant bit. 
We therefore pair the $2^{L+1}$ matrices 
$\mat{Z}^{m+1}(k_{m+1}, :, :) \in \mathbb{R}^{r \times r}$ over 
bit position $m$ in $k_{m+1}$: for each pair $(k^{(0)}, k^{(1)})$ 
differing only at bit $m$ (with $k^{(0)}$ having bit $m$ equal to $0$ and 
$k^{(1)}$ having bit $m$ equal to $1$), we concatenate along the left-rank 
axis to form a $(2r \times r)$ matrix and compute its reduced QR 
factorization. The two halves of $\mat{Q}$ are assigned back to indices 
$k^{(0)}$ and $k^{(1)}$, yielding $2^{L+1}$ matrices with pairwise 
orthonormal columns, while the remainder 
$\mat{R} \in \mathbb{R}^{r \times r}$ must be absorbed into 
the adjacent factor $\tsr{Z}^{m+2}$.

\paragraph{Absorption of $\mat{R}$}
Each remainder $\mat{R}$ is left-multiplied into the matrices of 
$\tsr{Z}^{m+2}$ that share the same remaining indices. Comparing the 
index sets of $\tsr{Z}^{m+1}$ and $\tsr{Z}^{m+2}$, the transition from 
level $m+1$ to $m+2$ drops $i_{L-m}$ and introduces $j_m$. The new index 
$j_m$ occupies bit position $0$ (the least significant bit) in 
$k_{m+2} = \psi^{-1}(i_0, \dots, i_{L-m-1}, j_0, \dots, j_m)$. Each 
remainder $\mat{R}$ from the QR step therefore multiplies both matrices 
of the pair $(k^{(0)}, k^{(1)})$ differing at bit $0$ in $k_{m+2}$:
\begin{align}
    \mat{Z}^{m+2}(k^{(b)}, :, :) \gets \mat{R} \, \mat{Z}^{m+2}(k^{(b)}, :, :), 
    \quad b \in \{0, 1\}.
\end{align}
When the target is the right outer factor $\mat{Z}^{L+2}$, each block 
$\mat{Z}^{L+2}_{k_{L+2}}$ of size $(r \times c)$ is instead updated as 
$\mat{Z}^{L+2}_{k_{L+2}} \gets \mat{R}_{k_{L+2}} \, \mat{Z}^{L+2}_{k_{L+2}}$.

The complete ADF algorithm for butterfly matrix completion is summarized in 
Algorithm~\ref{alg:butterfly-adf}. The gradient computation is the same as that used in~\cref{alg:ADAM_BF} and the optimal step size for each slice of the 
computed tensors are computed as suggested 
in~\cite{grasedyck2015variants}.

\begin{algorithm}[H]
\caption{\hl{ADF for Butterfly Completion}}
\label{alg:butterfly-adf}
\begin{algorithmic}[1]
\State \textbf{Input:} Observed matrix $\mat{T}_\Omega$ (or its tensorization $\tsr{T}_\Omega$), rank \( r \), initial guess $\tsr S^{i},i=1,\ldots,L+2$, max iterations \( t_{\max} \), and convergence threshold $\epsilon$.
\State \textbf{Output:} Updated butterfly cores $\tsr S^{i},i=1,\ldots,L+2$.
\State $\mathcal{Z} \gets \textsc{RightToLeftSweep}(\mathcal{Z}, L)$
\Comment{LQ sweep: concentrate weight in $\tsr{Z}^1$}
\State $\tsr{G} \gets \tsr{T}_\Omega - 
\textsc{Reconstruct}_\Omega(\mathcal{Z})$
\Comment{Initial residual on $\Omega$}
\For{$\text{iter} = 1, \dots, t_{\max}$}
    \For{$l = 1, 2, \dots, L+2$}
        \State $\tsr{N} \gets \textsc{Gradient}(\tsr{G}, \Omega, 
        \tsr{Z}, l, L)$\label{line:adf:gradient}
        \Comment{ As in~\cref{alg:ADAM_BF}}
        \State $\tsr{Z} \gets \textsc{Reconstruct}_\Omega(
        \tsr{Z}|_{\tsr{Z}^l \leftarrow \tsr{N}})$\label{line:adf:reconstruct}
        \State $\alpha_s \gets \|\tsr{N}_{(s)}\|_F^2 \,/\, \|\vcr{z}\|_2^2$
        \Comment{Overrelaxation for each slice as 
        in~\cite{grasedyck2015variants}}
        \State $\tsr{Z}^l_s \gets \tsr{Z}^l_{s} + \alpha_s \, \tsr{N}_s$ \Comment{Slice-wise update}
        \State $\tsr{G}_s \gets \tsr{G}_s - \alpha_s \, \tsr{Z}_s$ \Comment{Slice-wise update}
        \If{$l \leq L + 1$} \Comment{Batched QR and absorb into next level}
            \State $\tsr{Z}^l, \mat{R} \gets 
            \textsc{QR}(\tsr{Z}^l, L, l)$
            \State $\tsr{Z}^{l+1} \gets 
            \textsc{Absorb}(\mat{R}, \tsr{Z}^{l+1}, L, l)$
        \EndIf
    \EndFor
    \If{$\|\vcr{g}\|_2 / \|\tsr{T}_\Omega\|_F < \epsilon$}
        \State \textbf{break}
    \EndIf
\EndFor
\end{algorithmic}
\end{algorithm}

\bibliographystyle{plainurl}
\bibliography{main}

\end{document}